\documentclass[a4paper,10pt,reqno]{article}

\usepackage[utf8]{inputenc}
\usepackage[T1]{fontenc}
\usepackage{lmodern}
\usepackage[english]{babel}
\usepackage{microtype}

\usepackage{amsmath,amssymb,amsfonts,amsthm}
\usepackage{mathtools,accents}
\usepackage{mathrsfs}
\usepackage{aliascnt}
\usepackage{braket}
\usepackage{bm}

\usepackage[a4paper,margin=3cm]{geometry}
\usepackage[citecolor=blue,colorlinks]{hyperref}

\usepackage{enumerate}
\usepackage{xcolor}
\makeatletter
\g@addto@macro\@floatboxreset\centering
\makeatother

\makeatletter
\def\newaliasedtheorem#1[#2]#3{
	\newaliascnt{#1@alt}{#2}
	\newtheorem{#1}[#1@alt]{#3}
	\expandafter\newcommand\csname #1@altname\endcsname{#3}
}
\makeatother

\numberwithin{equation}{section}

\newtheoremstyle{slanted}{\topsep}{\topsep}{\slshape}{}{\bfseries}{.}{.5em}{}

\theoremstyle{plain}
\newtheorem{theorem}{Theorem}[section]
\newaliasedtheorem{proposition}[theorem]{Proposition}
\newaliasedtheorem{lemma}[theorem]{Lemma}
\newaliasedtheorem{corollary}[theorem]{Corollary}
\newaliasedtheorem{counterexample}[theorem]{Counterexample}

\theoremstyle{definition}
\newaliasedtheorem{definition}[theorem]{Definition}
\newaliasedtheorem{question}[theorem]{Question}
\newaliasedtheorem{conjecture}[theorem]{Conjecture}

\theoremstyle{remark}
\newaliasedtheorem{remark}[theorem]{Remark}
\newaliasedtheorem{example}[theorem]{Example}

\newcommand{\setN}{\mathbb{N}}

\newcommand{\setR}{\mathbb{R}}

\newcommand{\eps}{\varepsilon}

\let\altphi\phi
\let\phi\varphi
\let\varphi\altphi
\let\altphi\undefined

\newcommand{\abs}[1]{\left\lvert#1\right\rvert}
\newcommand{\norm}[1]{\left\lVert#1\right\rVert}

\let\div\undefined
\DeclareMathOperator{\div}{div}

\newcommand{\di}{\mathop{}\!\mathrm{d}}

\DeclareMathOperator{\supp}{supp}

\newcommand{\Ch}{{\sf Ch}}

\DeclareMathOperator{\Lip}{Lip}

\DeclareMathOperator{\lip}{lip} % %for the slope

\newcommand{\leb}{\mathscr{L}}

\newcommand{\Prob}{\mathscr{P}}
\newcommand{\Borel}{\mathscr{B}}

\newcommand{\XX}{{\boldsymbol{X}}}

\newcommand{\dist}{\mathsf{d}}

\newcommand{\meas}{\mathfrak{m}}

\newcommand{\Test}{\rm Test}
\newcommand{\TestV}{\rm TestV}
\DeclareMathOperator{\CD}{CD}
\DeclareMathOperator{\RCD}{RCD}
\newfont{\tmpf}{cmsy10 scaled 2500}

\def\Xint#1{\mathchoice
	{\XXint\displaystyle\textstyle{#1}}%
	{\XXint\textstyle\scriptstyle{#1}}%
	{\XXint\scriptstyle\scriptscriptstyle{#1}}%
	{\XXint\scriptscriptstyle\scriptscriptstyle{#1}}%
	\!\int}
\def\XXint#1#2#3{{\setbox0=\hbox{$#1{#2#3}{\int}$ }
		\vcenter{\hbox{$#2#3$ }}\kern-.6\wd0}}

\def\dashint{\Xint-}

\begin{document}
	
	\title{Regularity of Lagrangian flows over $\RCD^*(K,N)$ spaces}
	\author{Elia Bru\'e, Daniele Semola \thanks{Scuola Normale Superiore, \url{elia.bru\'e@sns.it, daniele.semola@sns.it}}} 
	\maketitle
		
	\begin{abstract}
		The aim of this note is to provide regularity results for Regular Lagrangian flows of Sobolev vector fields over compact metric measure spaces verifying the Riemannian curvature dimension condition. 
		
		We first prove, borrowing some ideas already present in the literature, that flows generated by vector fields with bounded symmetric derivative are Lipschitz, providing the natural extension of the standard Cauchy-Lipschitz theorem to this setting. Then we prove a Lusin-type regularity result in the Sobolev case (under the additional assumption that the m.m.s. is Ahlfors regular) therefore extending the already known Euclidean result.
		  
	\end{abstract}

	\tableofcontents
	\section*{Introduction}
	The theory of metric measure spaces with Riemannian Ricci curvature bounded from below and dimension bounded from above ($\RCD^*(K,N)$ metric measure spaces for short), although being very recent, is a very rapidly increasing research area with several contributions that, apart from their own theoretical interest, have often given new insights in the understanding of more classical questions of analysis and Riemannian geometry. 
	
	The introduction of the notion of metric measure spaces with Ricci curvature bounded from below and dimension bounded from above ($\CD(K,N)$ m.m.s. for short) dates back to the seminal and independent works of Lott-Villani \cite{LottVillani} and Sturm \cite{Sturm06a}, \cite{Sturm06b}. Crucial properties of this theory (which is therein formulated in terms of convexity-type properties of suitable energies over the Wasserstein space) are the compatibility with the case of smooth Riemannian manifolds and the stability w.r.t. suitable notions of convergence of metric measure spaces.
	
	Many geometrical and analytical properties have been proven for $\CD(K,N)$ metric measure spaces (see for instance \cite{Villani09}). However this class turns to be too large for some purposes, since it includes for instance smooth Finsler manifolds.
	In order to single out spaces with a Riemannian like behaviour from the above introduced broader class, in \cite{AmbrosioGigliSavare14} the authors proposed a notion of m.m.s. with Riemannian Ricci curvature bounded from below,
	adding to the $\CD$ condition the requirement of linearity of the heat flow (which is the gradient flow of the so-called Cheeger energy). Later on the theory was adapted in \cite{Gigli15}, \cite{AmbrosioMondinoSavare15} and \cite{ErbarKuwadaSturm15} to the dimensional case, with the introduction of the $\RCD(K,N)$ condition
	\footnote{We will work with the slightly modified $\RCD^*(K,N)$ condition in this paper. We avoid any further comment about the differences w.r.t. the $\RCD(K,N)$ condition in this introduction but we remark that very recently Cavalletti and Milman proved their equivalence in \cite{CavallettiMilman16}}.

	\paragraph*{}
	
	This paper deals with the regularity of flows of vector fields over $\RCD^*(K,N)$ metric measure spaces. To better introduce the reader to the notions that will be considered in the rest of the paper we briefly recall the Euclidean side of the story, which has been considered from much more time in the literature (but still deserves challenging open problems and questions). 
	
	\paragraph*{}
	
	In the Euclidean setting the Cauchy-Lipschitz theory grants existence, uniqueness and Lipschitz regularity for flows of Lipschitz vector fields. It is well known instead that lowering the regularity assumptions on the vector field might lead to non-uniqueness for integral curves, moreover, if one considers vector fields that are not defined everywhere but only Lebesgue-almost everywhere there is also need to introduce a notion of flow more general w.r.t. the one adopted in the smooth case.
	
	Motivated by the study of some PDEs in kinetic theory and fluid mechanics, Di Perna and Lions introduced in \cite{lions} a suitable notion of flow of Sobolev vector field and studied the associated existence and uniqueness problem. Later on their theory was revisited and extended to the case of vector fields with BV spatial regularity by Ambrosio in \cite{Ambrosio04}, where the notion of Regular Lagrangian Flow was introduced as a good global selection of integral curves of the vector field. Moreover Crippa and De Lellis in \cite{CrippaDeLellis08} were able to prove a mild regularity result for Regular Lagrangian Flows of Sobolev vector fields, namely that (locally) they are Lipschitz if we neglect a subset of the domain whose measure can be made arbitrary small (where, of course, we have to pay the price that the Lipschitz constant becomes arbitrary large). Such a result, known in the literature as Lusin-type regularity, holds true for instance for real valued Sobolev functions (and it is already known to be true also when the domain is a sufficiently regular metric measure space, see \cite{AmbrosioColomboDiMarino15}).
	
	\paragraph*{}
	
	Over an arbitrary metric measure space $(X,\dist,\meas)$ vector fields can be defined both as derivations over an algebra of scalar functions (which is the interpretation adopted in \cite{AmbrosioTrevisan14}) and as sections of the tangent modulus (see \cite{Gigli14} for the latter viewpoint and for the equivalence with the first one).
	Moreover, restricting the analysis to more regular metric measure spaces (such as $\RCD(K,\infty)$ metric measure spaces) where a second order differential calculus can be developed, one can introduce also reasonable notions of Sobolev vector fields (see again \cite{AmbrosioTrevisan14} and \cite{Gigli14} for the definitions of the spaces of vector fields with symmetric derivative in $L^2$ and of Sobolev vector fields, respectively).
	
	A remark concerning the discussion above is in order. On metric measure spaces we do not have a priori at our disposal a notion of Lipschitz vector field (and also in the case of smooth Riemannian manifolds this notion is less natural and more subtle, since it requires parallel transport to compare tangent vectors at different points); in addition we do not have a notion of tangent vector at a given point. With this said, when trying to develop a theory of flows of vector fields in this very abstract setting, it is more natural to look at the generalized theory of Regular Lagrangian Flows than at the Cauchy-Lipschitz theory. In \cite{AmbrosioTrevisan14} Ambrosio and Trevisan were able to prove that this theory was the right one in order to get existence and uniqueness of flows of vector fields with symmetric covariant derivative in $L^2$ on a large class of metric measure spaces including that one of $\RCD(K,\infty)$ spaces. 
	
	After having established such a result one might wonder if the flow maps have some further regularity property. At a speculative level, proving such a result for a certain class of metric measure spaces could give new insights about their geometry and their regularity. We remark that this very recent theory has already been useful in some applications (see \cite{GigliRigoni17}, \cite{Han17} and \cite{GigliKetterKuwadaOhta17}).   
	
	The main result of this note is \autoref{thm: sobolev vectorfield, regularity} below, where we extend the Lusin-type regularity result by Crippa-De Lellis from the Euclidean setting to that one of compact $\RCD^*(K,N)$ metric measure spaces $n$-Ahlfors regular for some $1<n\le N<+\infty$.  
	We remark that, up to our knowledge, this is also the first intrinsic proof of the regularity result over smooth compact Riemannian manifolds, since it avoids any use of local charts.
 % %aggiunto questo commento sulla questione delle varietà riemanniane.
	
	\paragraph*{}
	
	This paper is organized as follows: in the preliminary \autoref{sec:preliminaries}  we introduce the main notations and collect the basic results of the theory of $\RCD$ metric measure spaces and Regular Lagrangian flows that are needed in the rest of the work. Then, in \autoref{sec:reglipcase}, following some ideas already present in \cite{Han17} and \cite{Sturm14}, we provide a full Lipschitz regularity result for flows of vector fields with bounded symmetric derivative. Finally, in \autoref{sec:Lusinregularity}, which is the core of this note, we prove the sought Lusin-type regularity results.

	{\bf Acknowledgements}
	The authors would like to thank Luigi Ambrosio for suggesting them the study of this problem and for kind and numerous comments and suggestions. They acknowledge the support of the PRIN2015 MIUR Project “Calcolo delle Variazioni”. 
	
	\section{Preliminaries}\label{sec:preliminaries}
	\subsection{RCD metric measure spaces}

	Throughout this note by metric measure space (m.m.s. for short) we mean a triple $(X,\dist,\meas)$ where $(X,\dist)$ is a complete and separable metric space and $\meas$ is a probability measure defined on the Borel $\sigma$-algebra of $(X,\dist)$.
	
	We shall adopt the standard metric notation: we will indicate by $B(x,r)$ the open ball of radius $r$ centred at $x\in X$, by $\Lip(X,\dist)$ the space of Lipschitz functions over $(X,\dist)$, by $\Lip f$ the Lipschitz constant of $f\in \Lip(X,\dist)$. Moreover we introduce the notation
	\begin{equation*}
	\lip f(x):=\limsup_{y\to x}\frac{\abs{f(x)-f(y)}}{\dist(x,y)},
	\end{equation*}
	for the so-called slope of a function $f:X\to\setR$.

	We will denote by $L^p(X,\meas)=L^p(X)=L^p$ the spaces of Borel $p$-integrable functions over $(X,\meas)$ for any $1\le p\le+\infty$ and by $L^0(X,\meas)$ the space of $\meas$-measurable functions over $X$.

	Unless otherwise stated from now on we assume $(X,\dist,\meas)$ to be a compact $\RCD^*(K,N)$ metric measure space for some $K\in\setR$ (lower bound on the Ricci curvature) and $1\le N<+\infty$ (upper bound on the dimension).	Let us assume without loss of generality that $\meas$ is fully supported on $X$, this assumption is justified by the fact that if $(X,\dist,\meas)$ is $\RCD^*(K,N)$ then so is $(\supp\meas,\dist,\meas)$. 
	
	We remark that the notion of $\RCD^*(K,N)$ m.m.s. was introduced and firstly studied in \cite{Gigli15}, 
	%TODO: va formulato in modo migliore, in effetti in quel lavoro non c'è una "definizione" in cui Gigli dice "uno spazio RCD(K,N) è..." però varie volte accenna al fatto che questa potrebbe essere la condizione naturale per molte cose e di fatto ne prova delle conseguenze, stesso discorso per l'articolo sullo splitting, quindi secondo me va citato \cite{Gigli15}. Fra l'avevo messo come primo perché lo avevo trovato come prima citazione in qualche lavoro (forse di Gigli stesso però).
	\cite{AmbrosioMondinoSavare15} and \cite{ErbarKuwadaSturm15}), while the introduction of the $\RCD(K,\infty)$ condition dates back to the work \cite{AmbrosioGigliSavare14}. We just recall that those spaces can be introduced and studied both from an Eulerian point of view (based on the so-called $\Gamma$-calculus) and from a Lagrangian point of view (based on optimal transportation techniques). 
	
	Below we briefly describe the main analytic and geometric properties of $\RCD^*(K,N)$ metric measure spaces that will play a role in our work.

	As a first geometric property we recall that $\RCD^*(K,N)$ metric measure spaces satisfy the Bishop-Gromov inequality (which holds true more generally for any $\CD^*(K,N)$ m.m.s., see \cite{LottVillani}). Together with the compactness assumption, the Bishop-Gromov inequality implies that $(X,\dist,\meas)$ is doubling, that is there exists $c_D>0$ such that
	\begin{equation*}
	\meas(B(x,2r))\le c_D\meas(B(x,r))
	\end{equation*} 
	for any $x\in X$ and for any $r>0$.
	
	For any $f\in L^1(X,\meas)$ we will denote by $Mf$ the Hardy-Littlewood maximal function of $f$, which is defined by
	\begin{equation}\label{Maximal function}
	Mf(x):=\sup_{r>0}\dashint_{B(x,r)}\abs{f(z)}\di\meas(z),
	\end{equation} 
	where
	\begin{equation*}
	\dashint_{B(x,r)}f(z)\di\meas(z):=\frac{1}{\meas(B(x,r))}\int_{B(x,r)}f(z)\di\meas(z).
	\end{equation*} 
	We recall that, since $(X,\dist,\meas)$ is a doubling m.m.s., the maximal operator $M$ is bounded from $L^p(X,\meas)$ into itself for any $1<p\le+\infty$.
	
	We go on with a brief discussion about Sobolev functions and vector fields over $(X,\dist,\meas)$ referring to \cite{AmbrosioGigliSavare13}, \cite{AmbrosioGigliSavare14} and \cite{Gigli14} for a more detailed discussion about this topic. 
	
	\begin{definition}\label{def:Sobolevspace}
		For any $1<p<+\infty$ we define $W^{1,p}(X,\dist,\meas)(=W^{1,p}(X))$ to be the space of those $f\in L^{p}(X,\meas)$ such that there exists a sequence $f_n\to f$ in $L^p(X,\meas)$ with $f_n\in\Lip(X,\dist)$ for any $n\in\setN$ and $\sup_n\norm{\lip f_n}_{L^p}<+\infty$.
	\end{definition}
	
	The definition of Sobolev space is strongly related to the introduction of the Cheeger energy $\Ch_p:L^{p}(X,\meas)\to[0,+\infty]$ which is defined by
	\begin{equation}\label{eq:Cheeger}
	\Ch_p(f):=\inf\left\lbrace \liminf_{n\to\infty}\int \left(\lip f_n\right)^p\di\meas: f_n\to f\text{ in }L^p, f_n\in\Lip(X,\dist)\right\rbrace 
	\end{equation}
	and turns out to be a convex and lower semicontinuous functional from $L^p(X,\meas)$ to $[0,+\infty]$ whose finiteness domain coincides with $W^{1,p}(X,\dist,\meas)$.
	
	By looking at the optimal approximating sequence in \eqref{eq:Cheeger} one can identify a distinguished object, called minimal relaxed gradient and denoted by $\abs{\nabla f}_p$, which provides the integral representation
	\begin{equation*}
	\Ch_p(f)=\int_X\abs{\nabla f}_p^p\di\meas,
	\end{equation*}    
	for any $f\in  W^{1,p}(X,\dist,\meas)$.
	As the notation suggests $\abs{\nabla f}_p$ depends a priori on the integrability exponent $p$.
	
	The space $W^{1,p}(X,\dist,\meas)$ is a Banach space when endowed with the norm $\norm{f}_{W^{1,p}}^p:=\norm{f}^p_{L^p}+\Ch_p(f)$, moreover it holds that the inequality $\abs{\nabla f}_p\le\lip f$ holds true $\meas$-a.e. on $X$ for any $f\in\Lip(X,\dist)$. %and Lipschitz functions are dense in $W^{1,p}$ (see \cite{AmbrosioGigliSavare13}). 

	We point out that to single out $\RCD^*(K,N)$ metric measure spaces form the broader class of $\CD^*(K,N)$ metric measure spaces one adds the request that $\Ch:=\Ch_2$ is a quadratic form on $L^2(X,\meas)$ to the curvature-dimension condition. In this way the space $W^{1,2}(X,\dist,\meas)$, which in general is only a Banach space, turns to be a Hilbert space.
	
	This global assumption has in turn strong consequences on the infinitesimal behaviour of the space (indeed any m.m.s such whose $W^{1,2}$ is a Hilbert space is called infinitesimally Hilbertian). In particular, in the smooth setting it allows to single out Riemannian manifolds in the class of Finsler manifolds.
	
	For any $f,g\in W^{1,2}(X,\dist,\meas)$ we define a function $\nabla f\cdot\nabla g\in L^1(X\,\meas)$ by
	\begin{equation*}
    	\nabla f\cdot\nabla g(x):=\frac{1}{4}\abs{\nabla (f+g)}^2(x)-\frac{1}{4}\abs{\nabla (f-g)}^2(x) \quad\text{for $\meas$-a.e. $x\in X$.}
	\end{equation*} 

	It is proved in \cite{GigliHan16} that under the $RCD(K,\infty)$ assumption the minimal relaxed gradient $\abs{\nabla f}_p$ does not depend on $p$, for this reason we will use the notation $\abs{\nabla f}$.

	In order to introduce the heat flow and its main properties we begin by recalling the notion of Laplacian.
	
	\begin{definition}\label{def:laplacian}
		The Laplacian $\Delta:D(\Delta)\to L^2(X,\meas)$ is a densely defined linear operator whose domain consist of all functions $f\in W^{1,2}(X,\dist,\meas)$ satisfying 
		\begin{equation*}
		\int hg\di\meas=-\int \nabla h\cdot\nabla f\di\meas \quad\forall h\in W^{1,2}(X,\dist,\meas)
		\end{equation*}
		for some $g\in L^2(X,\meas)$. The unique $g$ with this property is denoted by $\Delta f$.\footnote{The linearity of $\Delta$ follows from the quadraticity of $\Ch$.}
	\end{definition}
	
	On any compact $\RCD^*(K,N)$ m.m.s the operator $-\Delta$ is densely defined, self-adjoint and compact. We will denote by $(\lambda_i)_{i\in\setN}$ its spectrum (where the eigenvalues are counted with multiplicity and in increasing order, $\lambda_i\ge 0$ for any $i$ and $\lambda_i\to\infty$ as $i$ goes to infinity) and by $(u_i)_{i\in\setN}$ the associated eigenfunctions, normalized in such a way that $\norm{u_i}_{L^2}=1$ for any $i\in\setN$. Recall that $(u_i)_{i\in\setN}$ is an orthonormal basis of $L^2(X,\meas)$.
	Furthermore the sequence $(\lambda_i)_{i\in\setN}$ has more than linear growth at infinity. A standard reference for this result in the smooth framework is \cite[Chapter 10]{Li12} and the arguments therein presented can be adapted to the case of our interest.  
	
	The heat flow $P_t$ is defined as the $L^2(X,\meas)$-gradient flow of $\frac{1}{2}\Ch$, whose existence and uniqueness follow from the Komura-Brezis theory. It can equivalently be characterized by saying that for any $u\in L^2(X,\meas)$ the curve $t\mapsto P_tu\in L^2(X,\meas)$ is locally absolutely continuous in $(0,+\infty)$ and satisfies
	\begin{equation*}
	\frac{\di}{\di t}P_tu=\Delta P_tu \quad\text{for $\leb^1$-a.e. $t\in(0,+\infty)$}.
	\end{equation*}  
	
	Under our assumptions the heat flow provides a linear, continuous and self-adjoint contraction semigroup in $L^2(X,\meas)$. Moreover $P_t$ extends to a linear, continuous and mass preserving operator, still denoted by $P_t$, in all the $L^p$ spaces for $1\le p<+\infty$.  
	
	In \cite{AmbrosioGigliSavare14} it is proved that for $\RCD(K,\infty)$ metric measure spaces the dual semgiroup $\bar{P}_t:\Prob(X)\to\Prob(X)$ of $P_t$, defined by
	\begin{equation*}
	    \int_X f \di \bar{P}_t \mu := \int_X P_t f \di \mu\qquad\quad \forall \mu\in \Prob(X),\quad \forall f\in \Lip_b(X),
	\end{equation*}
 for $t>0$, maps probability measures into probability measures absolutely continuous w.r.t. $\meas$. Then, for any $t>0$, we can introduce the so called \textit{heat kernel} $p_t:X\times X\to[0,+\infty)$ by
	\begin{equation*}
	p_t(x,\cdot)\meas:=\bar{P}_t\delta_x.
	\end{equation*}  
	From now on for any $f\in L^{\infty}(X,\meas)$ we will denote by $P_tf$ the representative pointwise everywhere defined by
	\begin{equation*}
	P_tf(x)=\int_X f(y)p_t(x,y)\di\meas(y).
	\end{equation*}
	
	Since $\RCD^*(K,N)$ metric measure spaces are doubling, as we already remarked, and they satisfy a local Poincaré inequality (see \cite{Villani09}) the general theory of Dirichlet forms (see \cite{Sturm96}) grants that we can find a locally H\"older continuous representative of $p$ on $X\times X\times(0,+\infty)$. 
	
	We also recall from \cite{AmbrosioHondaTewodrose17} the spectral identity which provides an explicit expression for the heat kernel in terms of the eigenfunctions of the Laplacian, namely
	\begin{equation}\label{eq:spectralidentity}
	p_t(x,y)=\sum_{i=0}^{\infty}e^{-\lambda_it}u_i(x)u_i(y),
	\end{equation}
	for any $t>0$, where, by choosing the H\"older continuous representative of $u_i$, whose H\"older norm growths linearly with $\lambda_i$, one obtains the H\"older continuous representative of $p_t$ (taking into account that, as we already remarked, the sequence of eigenvalues growths at least linearly at infinity).  
	
	Moreover in \cite{JangLiZhang} the following finer properties of the heat kernel have been proven: there exist constants $C_1\ge 1$, $C_2>0$ and $C_3\ge 0$ depending only on $K$ and $N$ such that 
	\begin{equation}\label{kernel estimates}
	\frac{1}{C_1}\frac{1}{\meas(B(x,\sqrt{t}))}\exp\left\lbrace-\frac{d^2(x,y)}{3t}-C_3t\right\rbrace
	\leq p_t(x,y)
	\leq C_1\frac{1}{\meas(B(x,\sqrt{t}))}\exp\left \lbrace-\frac{d^2(x,y)}{5t}+C_3t \right\rbrace,
	\end{equation}
	for any $x,y\in X$ and for any $t\in (0,+\infty)$ and 
	\begin{equation}\label{gradient of kernel estimate}
	|\nabla p_t(x,\cdot)|(y)\leq C_2\frac{1}{\meas(B(x,\sqrt{t}))\sqrt{t}}\exp\left\lbrace-\frac{d^2(x,y)}{5t}+C_3t\right\rbrace
	\qquad\quad \text{for}\ \meas\text{-a.e.}\ y\in Y,
	\end{equation}
	for any $t\in(0,+\infty)$ and for any $x\in X$.

	We go on by stating a few regularity properties of $\RCD^*(K,N)$ metric measure spaces (which hold true more generally for any $\RCD(K,\infty)$ m.m.s.) referring again to \cite{AmbrosioGigliSavare14} for a more detailed discussion and the proofs of these results.   
	
	First we have the \textit{Bakry-\'Emery contraction} estimate:
	\begin{equation}\label{eq:BakryEmery}
	\abs{\nabla P_t f}^2\le e^{-2Kt}P_t\abs{\nabla f}^2\quad \text{$\meas$-a.e.}
	\end{equation}
	for any $t>0$ and for any $f\in W^{1,2}(X,\dist,\meas)$.
	
	Another non trivial regularity property is the so-called \textit{$L^{\infty} - \Lip$ regularization} of the heat flow, that is for any $f\in L^{\infty}(X,\meas)$ we have that $P_tf\in\Lip(X)$ with
	\begin{equation}\label{eq:linftylipregularization}
	\sqrt{2I_{2K}(t)}\Lip(P_tf)\le\norm{f}_{L^\infty}\quad\text{for any $t>0$},
	\end{equation}
	where $I_L(t):=\int_0^te^{Lr}\di r$.
	
	Then we have the so-called \textit{Sobolev to Lipschitz property}: any $f\in W^{1,2}(X,\dist,\meas)$ with $\abs{\nabla f}\in L^{\infty}(X,\meas)$ admits a Lipschitz representative $\bar{f}$ such that $\Lip f\le \norm{\nabla f}_{\infty}$, that actually implies the equality $\Lip f= \norm{\nabla f}_{\infty}$ since, in general, $\Lip f\geq \norm{\lip f}_{\infty}$ and $\lip f\geq |\nabla f|$ $\meas$- a.e. on $X$.

	Following \cite{Gigli14} we  introduce the space of ``test'' functions $\Test(X,\dist,\meas)$ by 
	\begin{equation}\label{eq:test}
	\Test(X,\dist,\meas):=\{f\in D(\Delta)\cap L^{\infty}(X,\meas): \abs{\nabla f}\in L^{\infty}(X)\quad\text{and}\quad\Delta f\in W^{1,2}(X,\dist,\meas) \}
	\end{equation}
	and we remark that for any $g\in L^{\infty}(X)$ it holds that $P_tg\in \Test(X,\dist,\meas)$ for any $t>0$, thanks to \eqref{eq:BakryEmery}, \eqref{eq:linftylipregularization}, the fact that $P_t$ maps $L^2(X,\meas)$ in $D(\Delta)$ and the commutation
	\begin{equation*}
		\Delta	P_t f= P_t \Delta f\qquad \quad \forall f\in D(\Delta),
	\end{equation*}
	between $\Delta$ and $P_t$.
	
	We conclude this preliminary section with some finer regularity properties which hold true under the stronger assumption that the metric measure space is Ahlfors regular.
	
	\begin{definition}\label{def:Ahlforsregularity}
		We say that $(X,\dist,\meas)$ is $n$-Ahlfors regular for some $1\le n<+\infty$ if there exist constants $0<c_1\le c_2$ such that
		\begin{equation}\label{eq:Ahlforsregularity}
		c_1 r^n\le\meas(B(x,r))\le c_2r^n
		\end{equation}
		for any $0<r<D$ and for any $x\in X$, where we denoted by $D$ the diameter of $X$.
	\end{definition}

	\begin{remark}\label{remark: dist integrability}
		Let us observe that assumption \eqref{eq:Ahlforsregularity} grants integrability of certain powers of the distance function, namely for any $x\in X$ and for any $\alpha<n$ we have that $y\mapsto \dist(x,y)^{-\alpha}$ is $\meas$-integrable. Indeed by Cavalieri's formula we have that		
		\begin{align*}
		\int_X \frac{1}{\dist(x,y)^{\alpha}} \di\meas(x) 
		=& \int_0^{\infty} \meas(\set{y\ :\ \dist(x,y)^{-\alpha}> \lambda }) \di \lambda\\
		=& \int_0^{\infty} \meas(B(x,\lambda^{-1/\alpha})) \di \lambda\\
		\leq & D^{-\alpha}+c_1\int_{D^{-\alpha}}^{\infty} \lambda^{-n/\alpha} \di \lambda.
		\end{align*}
	\end{remark}
	
	In $\RCD^*(K,N)$ spaces it can be proved that eigenfunctions of the Laplacian have Lipschitz representatives. The result of the forthcoming \autoref{lemma:estimatesforeigenfunctions} provides also quantitative estimates on their Lipschitz norms, under the additional assumption that $\meas$ is Ahlfors regular.

	\begin{lemma}\label{lemma:estimatesforeigenfunctions} 
		Let $u_i$ be an eigenfunction of $-\Delta$ associated to the eigenvalue $\lambda_i$. Then $u_i$ has a Lipschitz representative. Moreover it holds
		\begin{equation}\label{eigenfunctions estimate}
		\norm{u_i}_{L^{\infty}}\leq \frac{C_1 e}{c_1} (C_3+\lambda_i)^{n/2},
		\quad \quad
		\norm{\nabla u_i}_{L^\infty}\leq \sqrt{(\lambda_i+ |K|)/2} \norm{u_i}_{L^{\infty}}.	
		\end{equation}
	\end{lemma}
	\begin{proof}
		Observe that by \eqref{eq:Ahlforsregularity} and \eqref{kernel estimates} we get
		$$
		p_t(x,y)\leq \frac{C_1 e^{C_3 t}}{c_1 t^{n/2}}\qquad\quad \forall x,y\in X,
		$$
		which yields the ultracontractivity property of the heat semigroup, namely
		\begin{equation}\label{ultracontractivity}
		\norm{P_t u}_{L^{\infty}} \leq \frac{C_1 e^{C_3 t}}{c_1 t^{n/2}}\norm{u}_{L^1}
		\qquad\quad \forall t>0,
		\end{equation}
		for any $u\in L^1(X,\meas)$.
		
		Observe that, since $-\Delta u_i=\lambda_iu_i$, it holds that $P_tu_i=e^{-\lambda_it}u_i$ for any $t>0$.

		An application of \eqref{ultracontractivity} with $t=1/(C_3+\lambda_i)$ yields now to the desired estimate
		\begin{equation*}
		\norm{u_i}_{\infty} \leq \frac{C_1 e}{c_1} (C_3+\lambda_i)^{n/2},
		\end{equation*}
		since $\norm{u_i}_{L^1}\le\norm{u_i}_{L^2}=1$.
		
		In order to prove the second estimate in \eqref{eigenfunctions estimate} we apply \eqref{eq:linftylipregularization} to get
		\begin{equation*}
		\norm{\nabla u_i}_{L^{\infty}}=e^{\lambda_it}\norm{\nabla P_tu_i}_{L^{\infty}}=e^{\lambda_it}\Lip( P_tu_i)\le\frac{e^{\lambda_it}}{\sqrt{2I_{2K}(t)}}\norm{u_i}_{L^{\infty}}.
		\end{equation*} 
		Observing that $I_L(s)\ge se^{-\abs{L}s}$ and choosing $t:=1/(\lambda_i+\abs{K})$ we obtain the desired conclusion. 
		
	\end{proof}

	\subsection{Regular Lagrangian flows}\label{sec:introflows}
	
	In this subsection we recall the notion of regular Lagrangian flow (RLF for short) firstly introduced in the Euclidean setting by Ambrosio in \cite{Ambrosio04} (inspired by the earlier work by Di Perna and Lions \cite{lions}). 
	
	The notion of regular Lagrangian flow was introduced to study ordinary differential equations associated to weakly differentiable vector fields. It is indeed well-known that, in general, it is not possible to define in a unique way a flow associated to a non Lipschitz vector field since the trajectories starting from a fixed point are often not unique.
	Roughly speaking, the RLF is a selection of trajectories that provides a very robust notion of flow.
	
	In order to define the concept of regular Lagrangian flow over $\RCD^*(K,N)$ spaces
	we introduce the notion of vector field through that one of derivation that has been adopted in \cite{AmbrosioTrevisan14}.
	
	\begin{definition}\label{def:derivation}
		We say that a linear functional $b:\Lip(X,\dist)\to L^0(X,\meas)$ is a derivation if it satisfies the Leibniz rule, that is
		\begin{equation}\label{eq:Leibnizrule}
		b(fg)=b(f)g+fb(g),
		\end{equation}
		for any $f,g\in\Lip(X,\dist)$.
		
		Given a derivation $b$ we will write $\abs{b}\in L^p$ if there exist some function $g\in L^p(X,\meas)$ such that
		\begin{equation}\label{eq:continuityofderivation}
		b(f)\le g\abs{\nabla f} \quad\text{$\meas$-a.e. on $X$,}
		\end{equation}
		for any $f\in \Lip(X,\dist)$ and we will denote by $\abs{b}$ the minimal (in the $\meas$-a.e. sense) $g$ with such property.
	\end{definition}
	
	We will also use the notation $b\cdot\nabla f$ in place of $b(f)$ in the rest of the paper.
	
	We remark that if a derivation $b$ is in $L^p$ then it can be extended in a unique way to a linear functional on $W^{1,q}(X,\dist,\meas)$ still satisfying \eqref{eq:continuityofderivation}, where $q$ is the dual exponent of $p$.
	Moreover, any $f\in W^{1,2}(X,\dist,\meas)$ defines in a canonical way a derivation $b_f$ in $L^2$ through the formula $b_f(g)=\nabla f\cdot\nabla g$, usually called the \textit{gradient derivation} associated to $f$.
	
	A notion of divergence can be introduced by integration by parts.
	
	\begin{definition}\label{def:divergence}
		Let $b$ be a derivation with $\abs{b}\in L^1(X,\meas)$. We say that $\div b\in L^p(X,\meas)$ if there exists $g\in L^p(X,\meas)$ such that
		\begin{equation*}
		\int_X b(f)\di\meas=-\int_Xgf\di\meas
		\end{equation*}
		for any $f\in \Lip(X,\dist)$. By a density argument it is easy to check that such a $g$ is unique (when it exists) and we will denote it by $\div b$.
	\end{definition}

	In the rest of the note we will write $b\in L^p(TX)$ to denote a derivation such that $\abs{b}\in L^p(X,\meas)$. We refer to \cite{Gigli14} for the introduction of the so-called tangent and cotangent moduli over an arbitrary metric measure space and for the identification results between derivations and elements of the tangent modulus which stand behind the use of this notation.

	We also introduce here the notion of time dependent vector field.
	
	\begin{definition}\label{time dependent vector field}
		Let us fix $T>0$, and $p\in [1,+\infty]$. We say that a map $b:[0,T]\to L^p(TM)$ is a time dependent vector field if for every $f\in W^{1,q}(X,\dist,\meas )$ the map
		$$
		(t,x)\mapsto b_t\cdot \nabla f(x),
		$$
		is measurable with respect to the product sigma-algebra $\leb^1 \otimes \Borel(X)$.
		We say that $b$ is bounded if
		$$
		\norm{b}_{L^{\infty}} :=\norm{|b|}_{L^{\infty}([0,T]\times X)}<\infty,
		$$
		and that $b\in L^1((0,T);L^p(X,\meas))$ if
		$$
		\int_0^T \norm{|b_s|}_{L^p} \di s<\infty.
		$$
	\end{definition}
	
	In the context of $\RCD^*(K,N)$ spaces the definition of Regular Lagrangian flow reads as follows (see \cite{AmbrosioTrevisan14} and \cite{AmbrosioTrevisan15}). 	
	
	\begin{definition}\label{def:Regularlagrangianflow}
		Let us fix a time dependent vector field $b_t$ (see \autoref{time dependent vector field}). We say that a map $\XX:[0,T]\times X\rightarrow X$ is a Regular Lagrangian flow associated to $b_t$ if the following conditions hold true:
		\begin{itemize}
			\item [1)] $\XX(0,x)=x$ and $X(\cdot,x)\in C([0,T];X)$ for every $x\in X$;
			\item [2)] there exists a positive constant $L$, called compressibility constant, such that
			$$
			\XX(t,\cdot)_{\#} \meas\leq L\meas,
			$$
			for every $t\in [0,T]$;
			
			\item [3)] for every $f\in \Test(X,\dist,\meas)$ the map $t\mapsto f(\XX(t,x))$ is absolutely continuous on $[0,T]$ and
			$$
			\frac{\di}{\di t} f(\XX(t,x))= b_t\cdot \nabla f(\XX(t,x)) \quad \quad \text{for a.e.}\ t\in (0,T),\quad \text{for}\ \meas\text{-a.e}\ x\in X;
			$$
			.
		\end{itemize}
	\end{definition}
	The selection of ``good'' trajectories is encoded in condition 2), which is added to ensure that the trajectories of the flow do not concentrate with respect to the measure $\meas$.
	
	In the definition we are assuming that $\XX$ is defined in every point $x\in X$. Actually the notion of RLF is stable under modification in a negligible set of initial conditions, but we prefer to work with a pointwise defined map in order to avoid technical issues.
	
	The theory of Regular Lagrangian flows in the context of metric measure spaces was developed by Ambrosio and Trevisan in \cite{AmbrosioTrevisan14}. The authors work with a very weak notion of symmetric derivative for a vector field.
	
	\begin{definition}\label{def: symder A-T}
		Let $b\in L^2(TX)$ with $\div b\in L^2(X,\meas)$. We say that $|\nabla_{\text{sym}}b|\in L^2(X,\meas)$ if there exists a constant $c\geq 0$ satisfying
		\begin{equation}\label{sym der A-T formula}
		\left|\int_X \nabla_{\text{sym}}b(\nabla f, \nabla g) \di \meas \right|
		\leq c \norm{|\nabla f|}_{L^4}\norm{|\nabla g|}_{L^4},
		\end{equation}	
		for every $f,g\in \Test(X,\di,\meas)$, where
		\begin{equation*}
		\int_X \nabla_{\text{sym}}b(\nabla f, \nabla g) \di \meas
		:= -\frac{1}{2}\int_X \left\lbrace b\cdot \nabla f \Delta g+b\cdot \nabla g \Delta f-(\div b)\nabla f\cdot \nabla g \right\rbrace \di \meas.
		\end{equation*}
		We let $\norm{\nabla_{\text{sym}}b}_{L^2}$ be the smallest admissible $c$ in \eqref{sym der A-T formula}.
	\end{definition}
	
	The results of \cite{AmbrosioTrevisan14} grant in particular existence and uniqueness of the RLF associated to a bounded vector field $b$, with symmetric derivative (in the sense of \autoref{def: symder A-T}) in $L^2$ and bounded divergence, in the context of $\RCD^*(K,N)$ spaces, which is the one we will be interested on in the rest of this note.
	
	We conclude this section recalling a deep relation between the Regular Lagrangian flow of $b_t$ and solutions of the continuity equation induced by $b_t$ (see \cite{AmbrosioTrevisan14}).
	
	\begin{proposition}\label{prop:continuityregularlag}
		Let $\XX_t$ be a Regular Lagrangian flow of $b_t$. Then, for any $\mu_0\in\Prob(X)$ absolutely continuous w.r.t. $\meas$ and with density bounded from above, defining $\mu_t:=(\XX_t)_{\#}\mu_0$, we have that $\mu_t$ solves the continuity equation
		\begin{equation}\label{eq:continuityeq}
		\frac{\di}{\di t}\mu_t+\div(b_t\mu_t)=0
		\end{equation}
		in the distributional sense, that is the function $t\mapsto \int\varphi\di\mu_t$ is in $W^{1,1}(0,T)$ and it holds
		\begin{equation}\label{eq:continuity}
		\frac{\di}{\di t}\int\varphi\di\mu_t=\int b_t\cdot\nabla\varphi\di\mu_t,
		\end{equation}
		for $\leb^1$-a.e. $t\in(0,T)$, for any $\varphi\in\Lip(X,\dist)$.
	\end{proposition}

 \begin{remark}\label{rm:independencefromthetest}
 We remark that it is possible to find a common $\leb^1$-negligible set $\mathcal{N}\subset(0,T)$ in such a way that \eqref{eq:continuity} is satisfied for any $t\in(0,T)\setminus\mathcal{N}$ for any $\varphi \in\Lip(X,\dist)$ (see \cite[Proposition 3.7]{GigliHan15a})

 \end{remark}

	\section{Regularity in the ``Lipschitz'' case}\label{sec:reglipcase}
	
	The aim of this section is to provide a full Lipschitz regularity result for flows of (possibly time dependent) regular vector fields with bounded symmetric covariant derivative (where the right notions of symmetric covariant derivative and ``regular'' are introduced in \autoref{def:Sobolevvectorfield} below).
	
	We have to remark that, while the regularity assumption seems to be not too restrictive in view of the possible applications of this result, the assumption that the symmetric covariant derivative is bounded is very restrictive and it could happen that, for a general $\RCD^*(K,N)$ metric measure space, there are no vector fields satisfying this constraint. Nevertheless, we find it interesting to present this result, both to better introduce the reader to the study of flows of vector fields over non smooth spaces, both since techniques very similar to the one we are going to present have already proven to be useful in the study of some rigidity problems such as in \cite{GigliRigoni17} and \cite{GigliKetterKuwadaOhta17}.

	The proof of the sought regularity result follows the strategy of the Lipschitz regularity of flows of Lipschitz vector fields in the Euclidean case, based on the differentiation of the distance between two flow lines of the vector field. To rule out the possible non-smoothness of the space we work at the level of curves of absolutely continuous measures (exploiting the result of \autoref{prop:continuityregularlag}) following some ideas taken from \cite{Sturm14} and the very recent \cite{Han17}.
	Let us remark, for sake of correctness, that the strategy we implement, based on the application of the second order differentiation formula along $W_2$-geodesics, was already suggested in \cite{Han17} (see in particular Remark 3.18 and Remark 3.19 therein).
	
	We assume the reader to be familiar with the basic notions and notations of optimal transportation (referring for instance to \cite{AmbrosioGigliSavare13} for their introduction).
	
	Below we state two preliminary results that play a key role in the proof of \autoref{thm:quantitativeLipschitzforCE}. For the moment we don't need to add any extra regularity assumption to the time dependent vector field $(b_t)_{t\in[0,T]}$ apart from measurability w.r.t. time.

	\begin{lemma}\label{lemma:derivativealongflow}
		Let $(\mu_t)_{t\in[0,T]}$ be a solution of the continuity equation \eqref{eq:continuityeq} such that $\mu_t\le C\meas$ for any $t\in[0,T]$ for some $C>0$ and fix any $\nu\in\Prob(X)$. Then it holds
		\begin{equation}\label{eq:derivativeW_2}
		\frac{\di}{\di t}\frac{1}{2}W_2^2(\mu_t,\nu)=\int b_t\cdot\nabla\varphi_t\di \mu_t\quad\text{for $\leb^1$-a.e. $t\in (0,T)$},
		\end{equation}
		where $\varphi_t$ is any optimal Kantorovich potential for the transport problem between $\mu_t$ and $\nu$.
		
	\end{lemma}
	
	\begin{proof}
    The result is stated and proved in a much more general framework in \cite[Proposition 3.11]{GigliHan15a}.
	\end{proof}

	\begin{corollary}\label{cor:derivativealongdoubleflow}
		Let $(\mu_t)_{t\in[0,T]}$ and $(\nu_t)_{t\in[0,T]}$ be solutions with uniformly bounded densities to the continuity equation induced by $b_t$ starting from $\mu_0$ and $\nu_0$ respectively as in \eqref{eq:continuity}. Then it holds that
		\begin{equation}\label{eq:jointdistance}
		\frac{\di}{\di t}\frac{1}{2}W_2^2(\mu_t,\nu_t)\le\int b_t\cdot\nabla\varphi_t\di\mu_t+\int b_t\cdot\nabla\psi_t\di \nu_t\quad\text{for $\leb^1$-a.e. $t\in (0,T)$,}
		\end{equation}
		where $(\varphi_t,\psi_t)$ is any couple of optimal Kantorovich potentials between $\mu_t$ and $\nu_t$. 
	\end{corollary}

	\begin{proof}
    The desired conclusion can be obtained with almost the same proof of \autoref{lemma:derivativealongflow} above. We report here a few more details for sake of completeness.
    
    The results of \cite{GigliHan15a} grant that the curves $(\mu_t)_{t\in[0,T]}$ and $(\nu_t)_{t\in[0,T]}$ are absolutely continuous with values in $(\Prob(X), W_2)$. It easily follows that the curve $t\mapsto W_2^{2}(\mu_t,\nu_t)$ is absolutely continuous, hence it is differentiable $\leb^1$-a.e. over $(0,T)$.
    It follows from \autoref{rm:independencefromthetest} and what we just observed that we can find a full $\leb^1$-measure set $\mathcal{C}\subset(0,T)$ such that, for any $t\in\mathcal{C}$, $s\mapsto W^{2}_2(\mu_s,\nu_s)$ is differentiable at $s=t$ and it holds
    \begin{equation}\label{eq:derivatives}
    \frac{\di}{\di s}\vert_{s=t}\int\varphi\di\mu_s=\int b_t\nabla\cdot\varphi\di\mu_t,\quad\quad\frac{\di}{\di s}\vert_{s=t}\int\psi\di\nu_s=\int b_t\cdot\nabla\psi\di\nu_t.
    \end{equation}
    
    Let now $(\varphi_t,\psi_t)$ be any couple of optimal Kantorovich potentials between $\mu_t$ and $\nu_t$. It follows by the duality results for the optimal transport problem that for any $h>0$ sufficiently small it holds
    \begin{equation}\label{eq:duality}
    \frac{1}{2}W_2^2(\mu_{t+h},\nu_{t+h})-\frac{1}{2}W_2^2(\mu_t,\nu_t)\ge\int\varphi_t\di\mu_{t+h}+\int\psi_t\di\nu_{t+h}-\int\varphi_t\di\mu_t-\int\psi_t\di\nu_t.
    \end{equation} 
    The desired conclusion follows from \eqref{eq:duality} dividing by $h$, taking the limit as $h\to 0$ at both sides and taking into account \eqref{eq:derivatives}.
	\end{proof}

	Before going on we introduce following \cite{Gigli14} the notion of Sobolev vector field with symmetric covariant derivative in $L^2$ over an arbitrary $\RCD(K,\infty)$ m.m.s. $(X,\dist,\meas)$ and the associated space $W^{1,2}_{C,s}(X,\dist,\meas)$.
	
	We refer to \cite{Gigli14} for the construction of the modulus $L^2(T^{\otimes 2}X)$ (starting from the tangent modulus $L^2(TX)$) and for the introduction of the broader space $W^{1,2}_C(X,\dist,\meas)$ of vector fields with full covariant derivative in $L^2$.
	
	\begin{definition}\label{def:Sobolevvectorfield}
		The Sobolev space $W^{1,2}_{C,s}(TX)\subset L^2(TX)$ is the space of all $b\in L^2(TX)$ for which there exists a tensor $S\in L^{2}(T^{\otimes 2}X)$ such that for any $h,g_1,g_2\in\Test(X,\dist,\meas)$ it holds
		\begin{equation}\label{eq:Sobvectfield}
		\int h S(\nabla g_1,\nabla g_2)\di\meas=\frac{1}{2}\int\left\lbrace -b(g_2)\div(h\nabla g_1)-b(g_1)\div(h\nabla g_2)+\div(hb)\nabla g_1\cdot\nabla g_2\right\rbrace \di\meas.
		\end{equation}
		In this case we shall call the tensor $S$ symmetric covariant derivative of $b$ and we will denote it by $\nabla_{\text{sym}} b$. 
		We endow the space $W^{1,2}_{C,s}(TX)$ with the norm $\norm{\cdot}_{W^{1,2}_{C,s}(TX)}$ defined by
		\begin{equation*}
		\norm{b}^2_{W^{1,2}_{C,s}(TX)}:=\norm{b}^2_{L^2(TX)}+\norm{\nabla_{\text{sym}} b}^2_{L^2(T^{\otimes 2}TX)}.
		\end{equation*}
	\end{definition}

	\begin{remark}
		It easily follows from the definition that the symmetric covariant derivative is actually symmetric.
		
		Moreover, for any $b\in W^{1,2}_C(TX)$ it holds that $b\in W^{1,2}_{C,s}(TX)$ and $\nabla_{\text{sym}}b$ is the symmetric part of $\nabla b$ (we refer to \cite[Section 3.4]{Gigli14} for the definition of the covariant derivative). 
		
	\end{remark}
	
	In order to compare the notion of Sobolev vector field introduced above with that one introduced in \autoref{def: symder A-T} we observe that the first one is easily seen to be stronger than the second one. Indeed it is sufficient to take $h=1$ in \autoref{eq:Sobvectfield} and to apply the Young inequality with exponents $2,4$ and $4$ to get the claimed conclusion.

	We define the space $\TestV(X,\dist,\meas)\subset L^2(TX)$ of test vector fields to be the set of linear combinations of the form
	\begin{equation*}
	\sum_{i=1}^ng_i\nabla  f_i,
	\end{equation*}
	where $g_i,f_i\in\Test(X,\dist,\meas)$ for any $i=1,\dots,n$.
	
	We recall that $\TestV(X)\subset W^{1,2}_C(X)$ and we will denote by $H^{1,2}_{C,s}(X)$ the closure of $\TestV(X)$ in $W^{1,2}_{C,s}(X)$ w.r.t. the $W^{1,2}_{C,s}$-norm.

	Below we quote from \cite{GigliTamanini17b} a useful differentiation formula for  $H^{1,2}_{C,s}$-vector fields along Wasserstein geodesics obtained by the authors as a corollary of the second order differentiation formula for $H^{2,2}$-Sobolev functions proven in \cite{GigliTamanini17}.

	% %Fra l'altro sembra che l'ipotesi di compattezza dello spazio sia stata rimossa

	\begin{theorem}\label{thm:diffforumlaforvectorfields}
		Let $(X,\dist,\meas)$ be a compact $\RCD^*(K,N)$ m.m.s. for some $1\le N<+\infty$. Let $(\eta_s)_{s\in[0,1]}$ be a $W_2$-geodesic connecting probability measures $\eta_0$ and $\eta_1$ absolutely continuous w.r.t. $\meas$ and with bounded densities and assume that $b\in H^{1,2}_{C,s}(X,\dist,\meas)$. Then one has that the curve
		\begin{equation*}
		s\mapsto \int b\cdot\nabla\varphi_s\di\eta_s
		\end{equation*}
		is $C^1$ on $[0,1]$, where $\varphi_s$ is any function such that for some $r\in[0,1]$ with $s\neq r$ it holds that $-(r-s)\varphi_s$ is an optimal Kantorovich potential from $\eta_s$ to $\eta_r$. Moreover it holds that
		\begin{equation*}
		\frac{\di}{\di s}\int b\cdot\nabla\varphi_s\di\eta_s=\int\nabla_{\text{sym}}b(\nabla\varphi_s,\nabla\varphi_s)\di\eta_s
		\end{equation*}
		for any $s\in[0,1]$.
	\end{theorem}
	
	\begin{remark}
		We remark that \autoref{thm:diffforumlaforvectorfields} is actually stated in \cite{GigliTamanini17b} only for vector fields in $H^{1,2}_{C}(X)$. However, since the strategy of the proof goes via approximation through elements of $\TestV(X)$ and the statement just involves the symmetric part of the covariant derivative, it easily extends to $H^{1,2}_{C,s}(X)$.
	\end{remark}

	We remark that it makes sense to say that an element of $L^2(T^{\otimes 2}X)$ belongs to $L^{\infty}(T^{\otimes 2}X)$ and to consider its $L^{\infty}$-norm.
	With this said, we will denote by 
	\begin{equation*}
	L:=\sup_{t\in (0,T)}\norm{\nabla_{\text{sym}}b_t}_{L^{\infty}}
	\end{equation*}
	and from now on to the basic assumptions of \autoref{sec:introflows} we add the assumption that $L<+\infty$.

	Below we state and prove the key result of this section, that will allow us to obtain both uniqueness and Lipschitz regularity for Regular Lagrangian Flows.
	
	\begin{theorem}\label{thm:quantitativeLipschitzforCE}
		Let $(X,\dist,\meas)$ be a compact $\RCD^*(K,N)$ m.m.s. and $(b_t)_{t\in[0,T]}$ be a time dependent vector field verifying the above discussed assumptions. Let $(\mu_t)_{t\in[0,T]}$ and $(\nu_t)_{t\in[0,T]}$ denote solutions of the continuity equation induced by $b_t$, absolutely continuous and with uniformly bounded densities. Then it holds that
		\begin{equation*}
		W_2(\mu_t,\nu_t)\le e^{Lt}W_2(\mu_0,\nu_0)
		\end{equation*} 
		for any $t\in[0,T]$.
	\end{theorem}

	\begin{proof}
		Applying first \autoref{cor:derivativealongdoubleflow} and then \autoref{thm:diffforumlaforvectorfields} we obtain that, for $\leb^1$-a.e. $t\in(0,T)$,
		\begin{align}\label{eq:estimateLip}
		\frac{\di}{\di t}\frac{1}{2}W_2^{2}(\mu_t,\nu_t)&\le\int b_t\cdot\nabla\mathcal{Q}_1\footnote{Here $\mathcal{Q}_t$ stands for the Hopf-Lax semigroup.}(-\varphi_t)\di\nu_t-\int b_\cdot\nabla-\varphi_t\di\mu_t\\
		&=\int_{0}^1\int\nabla_{sym}b_t(\nabla \varphi_t^s,\nabla\varphi_t^s)\di\eta_t^s\di s,
		\end{align}
		where we denoted by $\varphi_t$ an optimal Kantorovich potential from $\mu_t$ to $\nu_t$, $(\eta_t^{s})_{s\in[0,1]}$ the $W_2$-geodesic joining $\mu_t$ with $\nu_t$ and by $\varphi_t^s$ the intermediate time potential such that
		\begin{equation*}
		\frac{\di}{\di s}\eta_{t}^s+\div(\nabla\varphi_{t}^s\eta_{t}^s)=0.
		\end{equation*}
		Observe that we are in position to apply \autoref{thm:diffforumlaforvectorfields} since $\mu_t$ and $\nu_t$ are by assumption absolutely continuous with bounded densities for any $t\in[0,T]$.
		
		Recalling that, as a consequence of the metric Brenier theorem (see \cite[Proposition 3.5 ]{AmbrosioGigliSavare14}), 
		\begin{equation*}
		\int\abs{\nabla\varphi_t^s}^2\di\eta_t^s=W_2^2(\mu_t,\nu_t), 
		\end{equation*}
		for $\leb^1$-a.e. $s\in(0,1)$ and for any $t\in[0,T]$, we can conclude from \eqref{eq:estimateLip} that
		\begin{equation*}
		\frac{\di}{\di t}\frac{1}{2}W^{2}_{2}(\mu_t,\nu_t)\le \norm{\nabla_{\text{sym}}b_t}_{L^{\infty}}W_{2}^2(\mu_t,\nu_t)\le LW_2^2(\mu_t,\nu_t)
		\end{equation*}
		for a.e. $t\in(0,T)$.
		
		An application of Gromwall's lemma yields now the desired conclusion, namely that
		\begin{equation*}
		W_2^2(\mu_t,\nu_t)\le e^{2Lt}W^2_2(\mu_0,\nu_0)
		\end{equation*}
		for any $t\in[0,T]$. 
	\end{proof}
	
	As a corollary of \autoref{thm:quantitativeLipschitzforCE}, under our more restrictive assumptions about the regularity of the vector field, we can prove uniqueness (avoiding the theory of renormalized solutions) and Lipschitz regularity of Regular Lagrangian flows.
	
	\begin{theorem}\label{thm:uniquenessforLip}
		Let $(X,\dist,\meas)$ be a $\RCD^*(K,N)$ m.m.s. and let $(b_t)_{t\in[0,T]}$ satisfy the assumptions of \autoref{thm:quantitativeLipschitzforCE}. Then there exist a unique Regular Lagrangian flow $(\XX_t)_{t\in[0,T]}$ of $(b_t)_{t\in[0,T]}$.
	\end{theorem}
	
	\begin{proof}
		We do not give a complete proof of this statement. We just say here that \autoref{thm:quantitativeLipschitzforCE} gives uniqueness of solutions to the continuity equation induced by $(b_t)$ in the class of probability measures a.c. with respect to $\meas$ and with bounded density. Thus we are in position to proceed as in the proof of \cite[Theorem 7.7]{AmbrosioTrevisan15} to obtain uniqueness of the RLF.
	\end{proof}

	\begin{theorem}\label{thm:LipschitzforLip}
		Let $(X,\dist,\meas)$ and $(b_t)_{t\in[0,T]}$ be as before. Then for any $t\in[0,T]$ we can find a representative of the RLF $\XX_t$ satisfying the Lipschitz estimate
		\begin{equation*}
		\dist(\XX_t(x),\XX_t(y))\le e^{Lt}\dist(x,y),
		\end{equation*}
		for any $x,y\in X$.
	\end{theorem}
	
	\begin{proof}
		As for the proof of \autoref{thm:uniquenessforLip} above we do not give all the details.
		We just say here that the Lipschitz estimate for trajectories (which can be thought as solutions to the continuity equation starting from Dirac deltas) follows from \autoref{thm:quantitativeLipschitzforCE} from an approximation procedure whose details can be found for instance in the proof of \cite[Theorem 3.14]{Han17} or in \cite{Sturm14}. 
	\end{proof}

	\section{Regularity in the Sobolev case}\label{sec:Lusinregularity}
	
	In this section we prove a regularity property of regular Lagrangian flows associated to Sobolev vector fields in the context of (compact) Ahlfors regular $\RCD^*(K,N)$ metric measure spaces (see \autoref{def:Ahlforsregularity}). In order to better present the result and the main ideas of the proof we begin from the Euclidean setting, that is our starting point (even though non compactness requires some modification w.r.t. the strategy that we will adopt in the core of this section).
	
	In $(\setR^d, |\cdot|, \leb^d)$ the theory was developed by Crippa and De Lellis in \cite{CrippaDeLellis08} (implementing some ideas that were already present in \cite{AmbrosioLecumberryManiglia}), the main regularity result therein proved is the following one.
	
	\begin{theorem}\label{Th: Crippa-De Lellis Euclidean}
		Let $\XX_t$ be a regular Lagrangian flow associated to a time dependent vector field $b_t\in L^1((0,T);W^{1,p}(\setR^d;\setR^d))\cap L^{\infty}((0,T);L^{\infty}(\setR^n;\setR^n))$ with $p>1$, and fix $R>0$. For every $\varepsilon>0$ there exists a compact set $K\subset B_R$ such that $\leb^d(B_R\setminus K)< \varepsilon$ and 
		\begin{equation*}
		\Lip(\XX_t|_K)\leq\exp \left\lbrace \frac{C\left( 1+ \int_0^T \norm{\nabla b_t}_{L^p(B_{\tilde{R}})}\di t\right)}{\varepsilon^{1/p}}\right\rbrace ,
		\end{equation*}
		for any $t\in[0,T]$, where  $\tilde{R}:=R+T\norm{b}_{L^{\infty}}$ and $C$ depends only on $d,R,p$ and $L$.
	\end{theorem}
	
	The technique adopted in \cite{CrippaDeLellis08} is based on a priori estimates of the functionals
	\begin{equation*}
	Q_{t,r}(x):=\dashint_{B(x,r)} \log \left( 1+\left(\frac{ |\XX_t(x)-\XX_t(y)| }{r}\right) \right) \di y,
	\end{equation*}
	which represent a sort of non-convex, discrete Cheeger energies associated to $\XX$. 
	
	Any $L^p$ bound of the function $x\to \sup_{r>0} Q_{t,r}(x)$, depending only on the Sobolev norm of $b$, can be seen to imply a Lusin-type regularity property for $\XX_t$ similar to the one in \autoref{Th: Crippa-De Lellis Euclidean}.
	
	In order to find bounds for $Q_{t,r}$ one starts differentiating with respect to the time variable: 
	\begin{equation*}
	\frac{\di}{\di t} Q_{t,r}(x)
	=\dashint_{B(x,r)} \frac{\frac{\di}{\di t} |\XX_t(x)-\XX_t(y)| }{r+|\XX_t(x)-\XX_t(y)|} \di y
	\leq \dashint_{B(x,r)} \frac{|b_t(\XX_t(x))-b_t(\XX_t(y))|}{r+|\XX_t(x)-\XX_t(y)|} \di y.
	\end{equation*}
	To go on it suffices to recall the maximal estimate:
	\begin{equation}\label{lusin scalar}
	\frac{|b_t(\XX_t(x))-b_t(\XX_t(y))|}{|\XX_t(x)-\XX_t(y)|}\leq C(M |\nabla b|(\XX_t(x))+M|\nabla b| (\XX_t(y))).
	\end{equation}
	Now, using the assumption that $\XX$ has bounded compression and the $L^p$ integrability of the maximal operator for $p>1$, it is simple to find an $L^p$ bound of $x\to \sup_{r>0} Q_{t,r}(x)$ depending only on the Sobolev norm of $b$.

	Inequality \eqref{lusin scalar} in the Euclidean case follows applying the well-known Lusin-approximation property of scalar Sobolev functions with Lipschitz functions to the components of the vector field. Even though the scalar Lusin-approximation property is a very robust result (it holds true in every doubling metric measure space \cite{AmbrosioColomboDiMarino15} and in a rich class of non-doubling spaces \cite{AmbrosioBrueTrevisan17} ), it is non trivial to extend a similar property to vector fields out of the Euclidean setting.

	We are now ready to state the main result of this section. We refer to \autoref{thm: sobolev vectorfield, regularity} below for a more quantitative version of this statement.
	
	\begin{theorem} \label{thm: Lusin-type regularity of RLF}
		Let $(X,\dist, \meas)$ be a compact $\RCD^*(K,N)$ m.m.s. and  assume that it is $n$-Ahlfors regular for some $1<n\le N$ (see \autoref{def:Ahlforsregularity}).
		Let $(b_t)_{t\in[0,T]}$ be a bounded time dependent vector field with $\abs{\nabla_{\text{sym}} b_t}, \div b_t \in L^1((0,T); L^2(X,\meas))$. Let $\XX_t$ be a regular Lagrangian flow associated to $b_t$ with compressibility constant $L$. 
		Then for every $\epsilon>0$ there exists a Borel set $E\subset X$ such that
		$\meas(X\setminus E)<\epsilon$ and for every $x,y\in E$
		\begin{equation*}
		\dist(\XX_t(x),\XX_t(y))\leq C_{n,D}\exp\left\lbrace{C_{n,T} \frac{ L \int_0^t \norm{|\nabla_{\text{sym}}b_s|+|\div b_s|}_{L^2} \di s+1}{\sqrt{\epsilon}}}\right\rbrace \dist(x,y),
		\end{equation*}
		for every $t\in [0,T]$.
	\end{theorem}

	The Ahlfors regularity property, crucial in our proof, is a non trivial assumption. However the class of spaces we are able to treat is not poor, since it includes, for instance, Alexandrov spaces and non-collapsed $\RCD^*(K,N)$ metric measure spaces (see
	\cite{DePhilippisGigli17}).
	
	We conclude this preliminary discussion describing the main ideas in the proof of our result. Trying to perform the Crippa-De Lellis' 
	scheme the biggest difficulty to overcome comes from the study of the quantity 
	\begin{equation}\label{eq:derivativedistance}
	\frac{\di}{\di t}  \dist (\XX_t(x),\XX_t(y)). 
	%\leq | b\cdot \nabla \dist_{\XX_t(x)}(\XX_t(y))+b \cdot \nabla \dist_{\XX_t(y)}(\XX_t(x))|.
	\end{equation}	
	% %TODO: qua sopra noi dovremmo studiare la quantità a sinistra; se mettiamo il minore o uguale e poi studiamo la cosa a destra sembra che non riusciamo ad andare avanti perché abbiamo stimato male il LHS; visto che con la densità dell'algebra prodotto dovremmo riuscire a ottenere l'uguale potremmo pensare di metterlo nella prossima elaborazione. 
	% %Ho tolto il secondo pezzo della disuguaglianza per il motivo di cui parlo sopra
	
	Indeed in the metric setting it is not clear up to now how to obtain a useful estimate of the quantity
	\begin{equation}\label{difficult}
	b\cdot\nabla \dist_x(y)+b \cdot \nabla \dist_y(x)
	\end{equation}
	in terms of the covariant derivative of the vector field $b$ (a part from the case of bounded symmetric derivative which, however, seems to be too specific for the applications).
	
	Our strategy instead consists in considering a suitable power of the Green function $G(x,y)=G_x(y)$ instead of the distance function in \eqref{eq:derivativedistance} (we have been inspired by the survey \cite{ColdingMinicozzi12}).
	
	It is a well-known fact that on certain classes of Riemannian manifolds the Green function is equivalent to a negative power of the distance function; we extend this result to Ahlfors regular $\RCD^*(K,N)$ spaces.
	
	Instead of \eqref{difficult}, we need now to estimate
	\begin{equation}\label{less-difficult}
	b\cdot\nabla G_x(y)+b \cdot \nabla G_y(x)	;
	\end{equation}
	assuming for simplicity that $\div b=0$ and thanks to the fundamental property of the Green's function 
	\begin{equation*}
	\Delta G_x=\delta_x,
	\end{equation*}
	(actually $\Delta G_x=\delta_x-\meas$ in the case of compact manifolds)
	we formally compute
	\begin{align*}
	b\cdot\nabla G_x(y)+b \cdot \nabla G_y(x)
	& = \int_X [b(z)\cdot \nabla G_x(z) \Delta G_y(z)+b(z)\cdot \nabla G_y(z) \Delta G_x(z)] \di \meas(z)  \\
	& = -2 \int_X \nabla_{\text{sym}} b(\nabla G_x, \nabla G_y) \di \meas,
	\end{align*}
	that, with a little bit of work, provides a maximal estimate that plays the same role of \eqref{lusin scalar} in the Euclidean setting.
	
    As just said, throughout this section we make the additional assumption that $(X,\dist,\meas)$ is $n$-Ahlfors regular for some $1<n\le N$ (and it is still assumed to be a compact $\RCD^*(K,N)$ m.m.s. for some $1< N<+\infty$). 
	
	The main technical ingredients are developed in \autoref{subsec:Green}, where we prove that, for the class of spaces we are interested in, the Green function of the Laplace operator is comparable with a negative power of the distance function (extending a well-known result of Riemannian geometry, see \cite{Aubin98}). In \autoref{subsec:maximal} we turn the assumption on the Sobolev regularity of a vector field into a point-wise information obtaining a crucial maximal-type estimate. Through these two subsections we make the additional assumption that $n>2$, needed for technical reasons related to the different behaviour of the Green function in dimension two. Finally in \autoref{sec: n leq 2} we propose a short argument to extend the main result to the missing case $n=2$. We remark that, due to the results of \cite{MondinoNaber14}, the Ahlfors regularity assumption forces $n$ to be an integer between $1$ and $N$. Therefore the only remaining case would be that of $n=1$, that can be considered by iterating twice the procedure described in \autoref{sec: n leq 2}.
	% %Ho aggiunto un remark sul fatto che in realtà è solo $n=2$.
	
    % %in realtà è solo $n=2$, Mondino-Naber implica che se sei n-Ahlfors regolare allora n è intero (e minore o uguale a N)

	In order to let the notation be shorter we adopt the following convention: every positive constant that depends only on the ``structural'' coefficients of the space, i.e. $K,N,n,D, (\lambda_i)_{i\in \setN}, C_1,C_2,C_3,c_1,c_2$ and on universal numerical constants, will be denoted by $C$.

	\subsection{The Green function}\label{subsec:Green}
	Let us introduce now a key object for the rest of this note, namely the Green function
	\begin{equation}\label{eq:introGreen}
	G(x,y):=\int_0^{\infty} (p_t(x,y)-1) \di t, \quad\quad\quad \forall x,y\in X.
	\end{equation}
	In \autoref{prop: properties of G} below we prove that $G$ is well defined and we collect some important properties, extending to the case of our interest some known estimates in Riemannian geometry (see \cite{Aubin98} and \cite{Grygorian06}).
	
	Recall that we are assuming $n>2$.
	
	\begin{proposition}\label{prop: properties of G}
		The Green function $G$ in \eqref{eq:introGreen} is well defined and finite for every $x\neq y\in X$. For every $f\in \Test(X,\dist,\meas)$ it holds
		
		\begin{equation}\label{Green anction}
		\int_X G(x,y) \Delta f(y) \di\meas(y)= \int_X f \di\meas-f(x)
		\end{equation}
		for every $x\in X$.
		
		Moreover, $G$ is equivalent to the function $\dist(x,y)^{-n+2}$ up to a constant, i.e. there exist $A\geq 1$ and $\bar{A}>0$, depending only on $(X,\dist, \meas)$, such that
		\begin{equation}\label{G-d estimate 1}
		|G(x,y)|\leq \frac{A}{\dist(x,y)^{n-2}}, \quad \quad \quad \forall x,y\in X
		\end{equation}
		and 
		\begin{equation} \label{G-d estimate 2}
		G(x,y)\geq \frac{1}{A \dist(x,y)^{n-2}}-\bar{A} \quad \quad \quad \forall x,y\in X.
		\end{equation}
		
		Finally, setting $G_x(y):=G(x,y)$, there exists $C>0$ such that
		\begin{equation} \label{ nabla G estimate}
		\lip G_x(y)\leq \frac{C}{\dist(x,y)^{n-1}}\quad \quad \quad \forall x\neq y\in X,
		\end{equation}
		in particular $G_x, \lip G_x\in L^p(X,\meas)$ for every $p\in [1, n/(n-1))$.
		
	\end{proposition}
	
	\begin{proof}
		Let us prove that the integral $\int_0^{\infty} (p_t(x,y)-1) \di t$ is absolutely convergent.
		We assume for simplicity that the diameter $D$ of the space is equal to $1$.  Let us fix $x\neq y\in X$, using the estimates \eqref{kernel estimates} for the heat kernel, identity \eqref{eq:spectralidentity}, the Ahlfors regularity \eqref{eq:Ahlforsregularity} and \eqref{eigenfunctions estimate}, we have
		\begin{align*}
		\int_0^{\infty} |p_t(x,y)-1| \di t 
		&= \int_0^ 1 |p_t(x,y)-1| \di t+\int_1^{\infty}|p_t(x,y)-1| \di t\\
		&	\leq 1+\int_0^1 \frac{C_1}{\meas(B(x,\sqrt{t}))} e^{-\frac{\dist^2(x,y)}{5t}+C_3t} \di t+\int_1^{\infty}\sum_{i\geq 1} e^{-\lambda_i t} |u_i(x)||u_i(y)| \di t\\
		& \leq 1+\int_0^1 \frac{C}{t^{n/2}}e^{\frac{-\dist^2(x,y)}{5t}} \di t+\sum_{i\geq 1}\frac{Ce^{-\lambda_i}}{\lambda_i}(C_3+\lambda_i)^n.
		\end{align*}
		Observe now that the series $\sum_{i\geq 1}\frac{Ce^{-\lambda_i}}{\lambda_i}(C_3+\lambda_i)^n$ is convergent (since the eigenvalues have more than linear growth) and that 
		\begin{equation*}
		\int_0^1 \frac{C}{t^{n/2}}e^{-\frac{\dist^2(x,y)}{5t}} \di t\leq \frac{1}{\dist(x,y)^{n-2}}\int_0^{\infty} \frac{e^{-\frac{1}{5t}}}{t^{n/2}} \di t\leq \frac{C}{\dist(x,y)^{n-2}},
		\end{equation*}
		where in the last inequality the assumption that $n>2$ enters into play.
		All in all we have
		\begin{equation*}
		\left| \int_0^{\infty} (p_t(x,y)-1) \di t\right|\leq \frac{C}{\dist(x,y)^{n-2}}+C\leq \frac{C}{\dist(x,y)^{n-2}},
		\end{equation*}
		that provides the good definition of $G$ and \eqref{G-d estimate 1}.
		
		In order to prove \eqref{G-d estimate 2} observe that 
		\begin{equation*}
		\int_0^{\infty} (p_t(x,y)-1) \di t
		=\int_0^ 1 (p_t(x,y)-1) \di t+\int_1^{\infty}\sum_{i\geq 1} e^{-\lambda_i t} u_i(x)u_i(y) \di t.
		\end{equation*}
		Using again \eqref{kernel estimates} and \eqref{eq:Ahlforsregularity}, we obtain
		\begin{equation*}
		\int_0^ 1 (p_t(x,y)-1) \di t\geq \frac{C}{\dist(x,y)^{n-2}} -1.
		\end{equation*}
		Recalling that
		$$
		\left|\int_1^{\infty}\sum_{i\geq 1} e^{-\lambda_i t} u_i(x)u_i(y) \di t \right| 
		\leq \sum_{i\geq 1}\frac{Ce^{-\lambda_i}}{\lambda_i}(C_3+\lambda_i)^n
		<\infty
		$$ we conclude the proof of \eqref{G-d estimate 2}.
		
		Let us estimate now the slope of $G_x(\cdot)$ at $y\in X, \ y\neq x$. Let us fix a parameter $0<\epsilon< \dist(x,y)$, and a point $z\in B(y,\epsilon/2)$; observe that $\dist(x,z) >\epsilon/2$. We wish to estimate the incremental ratio
		\begin{equation*}
		\frac{|G_x(z)-G_x(y)|}{\dist(z,y)}\leq \int_0^1 \frac{|p_t(x,z)-p_t(x,y)|}{\dist(z,y)} \di t+ \int_1^{\infty} \sum_{i\geq 1} e^{-\lambda_i t}|u_i(x)|\frac{|u_i(z)-u_i(y)|}{\dist(y,z)} \di t=:I+II.
		\end{equation*}
		In order to estimate $I$ we observe that the slope of $p_t(x, \cdot)$ is bounded in $X\setminus B(x,\epsilon/2)$ uniformly in time and that a geodesic from $y$ to $z$ does not intersect $B(x,\epsilon/2)$. Thus, using the fact that the slope of a Lipschitz function is an upper gradient we obtain that the family 
		$\frac{|p_t(x,z)-p_t(x,y)|}{\dist(z,y)}$
		is uniformly bounded when $z\in B(y,\epsilon/2)$. By the dominated convergence theorem and \eqref{gradient of kernel estimate} we obtain
		\begin{align*}
		\limsup_{z\to y} \int_0^1 \frac{|p_t(x,z)-p_t(x,y)|}{\dist(z,y)} \di t 
		& \leq \int_0^1 |\nabla p_t(x,\cdot)|(y) \di t \\
		& \leq \int_0^1 \frac{C}{\meas(B(x,\sqrt{t}))\sqrt{t}} e^{-\frac{\dist^2(x,y)}{5t}} \di t \\
		&\leq \frac{C}{\dist(x,y)^{n-1}}.
		\end{align*}
		The estimate of $II$ is simple. Indeed from \eqref{eigenfunctions estimate} we obtain
		\begin{equation*}
		\limsup_{z\to y}\int_1^{\infty} \sum_{i\geq 1} e^{-\lambda_i t}|u_i(x)|\frac{|u_i(z)-u_i(y)|}{\dist(y,z)} \di t \leq 
		\int_1^{\infty} \sum_{i\geq 1} e^{-\lambda_i t}C(C_3+\lambda_i)^n(|K|+\lambda_i)^{1/2} \di t<+\infty.
		\end{equation*}
		Putting all together we conclude
		\begin{equation*}
		\lip G_x(y) \leq \frac{C}{\dist(x,y)^{n-1}}+C\leq \frac{C}{\dist(x,y)^{n-1}}.
		\end{equation*}
		By \autoref{remark: dist integrability} it easily follows that $G_x, \lip G_x \in L^p(X,\meas)$ for every $p\in [1,n/(n-1))$.
		
		Finally we prove \eqref{Green anction}. Let us fix $f\in \Test(X,\dist,\meas)$. We first
		observe that
		\begin{equation}\label{a1}
		\int_X G(\cdot,y) \Delta f(y) \di \meas(y)\in L^{\infty}(X,\meas),
		\end{equation}
		as a consequence of $\Delta f\in L^{\infty}(X,\meas)$ and \autoref{remark: dist integrability}.
%		\begin{align*}
%		\int_X |G(x,y)| |\Delta f(y)| \di \meas(y) \leq & C\sum_{k=-\infty}^0\int_{B(x,2^k)\setminus B(x,2^{k-1})}\frac{|\Delta f(y)|}{\dist(x,y)^{n-2}} \di \meas(y)\\
%		\leq & C\sum_{k=-\infty}^0\int_{B(x,2^k)\setminus B(x,2^{k-1})}\frac{|\Delta f(y)|}{2^{(n-2)(k-1)}} \di \meas(y)\\
%		\leq & C\sum_{k=-\infty}^0\frac{2^{kn}}{2^{(n-2)(k-1)}} \frac{1}{\meas(B(x, 2^k))} \int_{B(x,2^k)} |\Delta f(y)| \di \meas(y)\\
%		\leq & C M |\Delta f|(x),
%		\end{align*}
%		where $M$ is the maximal operator  introduced in \eqref{Maximal function}.
		Fix any $\phi\in L^2(X,\meas)$, applying Fubini's theorem we get
		\begin{align*}
		\int_X \phi(x)\int_X G(x,y) \Delta f(y)\di \meas(y)\di \meas(x)=&\int_0^\infty \int_X \phi(x)\int_X (p_t(x,y)-1)\Delta f(y) \di \meas(y) \di \meas(x) \di t\\
		=& \int_0^{\infty}\int_X  \phi(x)P_t\Delta f(x)\di \meas(x) \di t\\
		=& \int_0^{\infty} \frac{d}{dt}\int_X P_tf(x)\phi(x) \di \meas(x) \di t\\
		=&\int_X \left( \int_X f \di\meas-f(x)\right)\phi(x) \di \meas(x),
		\end{align*}    
	where all the integrals are well-defined thanks to \eqref{a1}.	 
	\end{proof}
	Let us introduce a ``regularized'' version $G^{\epsilon}$, $\epsilon>0$, of $G$ setting
	\begin{equation*}
	G^{\epsilon}(x,y):= \int_{\epsilon}^{\infty} (p_t(x,y)-1) \di t \quad \quad \quad \forall x,y\in X.
	\end{equation*}
	We will often write $G^{\epsilon}_x(y)=G^{\epsilon}(x,y)$. Observe that $G^{\epsilon}$ is well defined and finite for every $x,y\in X$. Estimates \eqref{G-d estimate 1} and \eqref{ nabla G estimate} still hold true for $G^{\epsilon}$, namely
	\begin{equation}\label{G eps estimate}
	|G^{\epsilon}(x,y)|\leq \frac{C}{\dist(x,y)^{n-2}},\qquad
	\lip G_x^{\epsilon}(y)\leq \frac{C}{\dist(x,y)^{n-1}}
	\end{equation}
	for every $x,y\in X$, and they can be proved with the the same strategy described above.
	In \autoref{action of G eps} below we state an important regularity property of $G^{\epsilon}_{x}$.
	\begin{lemma} \label{action of G eps}
		For every $x\in X$ it holds that $G^{\epsilon}_x \in \Test(X,\dist,\meas)$ and 
		\begin{equation}\label{a2}
		\Delta G^{\epsilon}_x (y)=1-p_{\epsilon}(x,y).
		\end{equation}
		for $\meas$-a.e. $y\in X$.
	\end{lemma}

	\begin{proof}
		Arguing as in the proof of \eqref{G-d estimate 1} we easily obtain $G_x^{\epsilon}\in L^{\infty}(X,\meas)$ and, with a simple application of Fubini-Tonelli theorem, we get 
		$$
		G_x^{\epsilon}=P_{\epsilon/2}G_x^{\epsilon/2}.
		$$
		Taking into account the regularizing properties of the heat flow that we remarked after \eqref{eq:test}, we obtain $G_x^{\epsilon}\in \Test(X,\dist,\meas)$.
		
		Identity \eqref{a2} follows arguing as in the proof of \eqref{Green anction}.
	\end{proof}
	Finally we observe that, for every $x\in X$, the family of functions $(G^{\epsilon}(x))_{\epsilon>0}$ is equibounded in $W^{1,p}(X,\dist,\meas)$ for some $p>1$, and strongly convergent as $\epsilon\to 0$ to $G_x$, details can be found in \autoref{cor: G sobolev} below.
	\begin{lemma}\label{cor: G sobolev}
		For every $x\in X$ and for every $p\in [1,\frac{n}{n-1})$ it holds that $G_x, G^{\epsilon}_x \in W^{1,p}(X,\dist,\meas)$ and 
		$$
		\lim_{\epsilon\to 0} G^{\epsilon}_x= G_x,\quad \quad  \quad \text{strongly in }\ W^{1,p}.
		$$
	\end{lemma}
	\begin{proof}
		From \eqref{G-d estimate 1}, \eqref{G eps estimate} and \autoref{remark: dist integrability} it immediately follows that
		$G_x, G^{\epsilon}_x \in L^p(X,\meas)$ for every $p\in [1,n/(n-2))$, with a bound on $L^p$ norms independent of $\epsilon$.
		Moreover $G_x^{\epsilon}$ is a Lipschitz function, thus
		$$
		|\nabla G_x^{\epsilon}|\le\lip G_x^{\epsilon}\leq 
		\frac{C}{\dist(x,y)^{n-1}},\qquad \text{$\meas$-a.e. on $X$}.
		$$
		Using \autoref{remark: dist integrability} again we conclude that
		$\sup_{\epsilon>0} \norm{G_x^{\epsilon}}_{W^{1,p}}< \infty$.
		It is simple to check that $G^{\epsilon}_x(y)\to G_x(y)$ for any $y\neq x$, when $\epsilon \to 0$ and by \eqref{G eps estimate} and the dominated convergence theorem we get
		\begin{equation}\label{a3}
		G_x^{\epsilon}\to G_x \qquad\quad \text{in}\ L^p\quad \text{for all } p\in [1,n/(n-2)).
		\end{equation}
		% %TODO: perché qua cambia l'esponente?
		Let us fix $p\in [1,\frac{n}{n-1})$.
		It is now obvious that $G_x\in W^{1,p}(X,\dist,\meas)$, since by the reflexivity of $W^{1,p}$ for $p>1$ (see \cite{AmbrosioColomboDiMarino15}) and \eqref{a3} we deduce that $(G_x^{\epsilon})_{\epsilon>0}$ weakly converges to $G_x$.
		It remains to prove that the convergence of $G_x^{\epsilon}\to G_x$ is actually strong in $W^{1,p}$. To this aim it is enough to show
		 $$
	   	 \lim_{\eps \to 0}  \int_X |\nabla G^{\eps}_x| \di \meas = \int_X |\nabla G_x| \di \meas,
	   	   $$
		since $G_x^{\eps}\to G_x$ weakly in $W^{1,p}$ and the space $W^{1,p}(X,\dist,\meas)$ is equivalent to a uniformly convex space (see \cite[Theorem 7.4]{AmbrosioColomboDiMarino15}). Using a $p$-version of \eqref{eq:BakryEmery} (see \cite[Proposition 3.1]{GigliHan16}) and the identity $P_{\eps}G_x=G^{\eps}_x$ (where the action of the semigroup is understood in $L^p(X,\meas)$) we have that
		\begin{equation*}
		     \int_X |\nabla G^{\eps}|^p(y)\di \meas (y)\leq e^{-pK\eps}\int_X |\nabla G|^p(y)\di \meas (y).
		\end{equation*}
		Therefore
		\begin{equation*}
			\limsup_{\eps\to 0}\int_X |\nabla G^{\eps}|^p(y)\di \meas (y)\leq \int_X |\nabla G|^p(y)\di \meas (y).
		\end{equation*} 
		Using the lower semicontinuity of the Sobolev norm with respect to the weak topology we conclude the proof.
	\end{proof}
%		for $p\in [1, \frac{n+2}{n+1})$. To this aim we show that $(G_x^{\epsilon})_{\epsilon>0}$ is a Cauchy sequence. Let us fix $\epsilon_1<\epsilon_2<1$, arguing as in proof of \autoref{prop: properties of G}, we obtain
%		$$
%		|\nabla (G^{\epsilon_2}_x-G^{\epsilon_1}_x)|(y)
%		\leq \lip(G^{\epsilon_2}_x-G^{\epsilon_1}_x)(y)
%		\leq C\int_{\epsilon_1}^{\epsilon_2}\frac{e^{-\frac{\dist^2(x,y)}{5t}+C_3t}}{\meas(B(x,t))\sqrt {t}} \di t,
%		$$
%		for $\meas$- a.e. $y\in X$. Thus 
%		\begin{align*}
%		\norm{|\nabla (G^{\epsilon_1}_x-G^{\epsilon_2}_x)|}_{L^p}^p 
%		&\leq C_2^pe^{C_3p}\int_X\left|\int_{\epsilon_1}^{\epsilon_2}\frac{e^{-\frac{\dist^2(x,y)}{5t}}}{\meas(B(x,\sqrt{t}))\sqrt{t}} \di t\right|^p \di \meas(y)\\
%		&\leq C_p|\epsilon_2-\epsilon_1|^{p-1} \int_{\epsilon_1}^{\epsilon_2} \frac{1}{t^{(n+1)p/2}}\int_X e^{-\frac{\dist^2(x,y)p}{5t}} \di \meas(y) \di t\\
%		&\leq C_p|\epsilon_2-\epsilon_1|^{p-1} \int_{\epsilon_1}^{\epsilon_2}\frac{\meas(B(x,\sqrt{5t/p}))}{t^{(n+1)p/2}} \di t\\
%		&\leq C_p |\epsilon_2-\epsilon_1|^{p-1}\int_{0}^{1} t^{-(n+1)p/2+n/2} \di t,
%		\end{align*}
%		that goes to zero when $\epsilon_1$ and $\epsilon_2$ go to zero for every $p\in [1, \frac{n+2}{n+1})$.

	\subsection{Maximal estimates for vector fields}\label{subsec:maximal}
	
	In this section we state and prove a maximal estimate which turns out to be crucial in the sequel. 
	\begin{proposition}\label{lemma: key estimate}
		There exists a positive constant $C=C(X,\dist,\meas)$ such that, for every $b\in L^2(TX)$ with $\div b,\nabla_{\text{sym}}b \in L^2(X,\meas)$, it holds
		\begin{equation*}
		|b\cdot \nabla G_x(y)+b \cdot \nabla G_y(x)|\leq
		\frac{C}{\dist(x,y)^{n-2}} [M(|\nabla_{\text{sym}}b|+|\div b|)(x)+M(|\nabla_{\text{sym}}b|+|\div b|)(y)],
		\end{equation*}
		for $\meas\times \meas$ a.e. $(x,y)\in X\times X$.
	\end{proposition}

	In some sense the result of \autoref{lemma: key estimate} could be seen as a quantitative Lusin-type approximation property for the vector valued case. Indeed it plays in our proof the role played by \eqref{lusin scalar} in the original Crippa-De Lellis' scheme (see \cite{CrippaDeLellis08}).
	
	The notion of symmetric covariant derivative we are adopting in \autoref{lemma: key estimate}, is the following one.	
	\begin{definition}\label{def: symder Elia-Dan}
		Take any $b\in L^2(TX)$ with $\div b\in L^2$. We say that $b$ has symmetric derivative in $L^2$ if there exists a non negative function $G\in L^2(X,\meas)$ such that
		\begin{equation}
		\frac{1}{2}\left|\int_X \left\lbrace b\cdot \nabla g \Delta f + b\cdot \nabla f \Delta g- \div b (\nabla f\cdot \nabla g)\right\rbrace  \di \meas\right|  \leq \int_X G |\nabla f||\nabla g| \di \meas,
		\end{equation}
		for every $f,g\in \text{Test(X)}$. We call $|\nabla_{\text{sym}}b|$ the $G\in L^2(X,\meas)$ with minimal norm\footnote{Up to now we do not know if $\abs{\nabla_{sym}b}$ has to be the minimal object also in the pointwise $\meas$-a.e. sense.}.	
	\end{definition}
	
	This definition appears as intermediate between the notion adopted by Ambrosio and Trevisan in \cite{AmbrosioTrevisan14} (see \autoref{def: symder A-T}) and the one proposed by Gigli in \cite{Gigli14} (see \autoref{def:Sobolevvectorfield}). Indeed it follows from the very definitions that if a vector field admits a symmetric covariant derivative in $L^2$ according to \autoref{def: symder Elia-Dan} then it also admits a covariant derivative in $L^2$ according to Ambrosio-Trevisan. On the other hand if a vector field belongs to $W^{1,2}_{C,s}$ (see \autoref{def:Sobolevvectorfield}), then it has symmetric covariant derivative in $L^2$ according to \autoref{def: symder Elia-Dan}.
	We chose to work with this intermediate notion of symmetric derivative since it is the assumption we really need for our purposes.
	
	We start with a technical lemma.
	\begin{lemma}\label{maximal estimate} There exists a positive constant $C=C(X,\di,\meas)$ such that for every non negative function $f\in L^1(X,\meas)$ it holds
		\begin{equation*}
		\int_X f(z)\frac{1}{\dist(x,z)^{n-1}}\frac{1}{\dist(y,z)^{n-1}} \di \meas(z) \leq \frac{C}{\dist(x,y)^{n-2}} \left( Mf(x)+ Mf(y)\right).
		\end{equation*}	 
	\end{lemma}

	\begin{proof}
		Set $r:=\dist(x,y)/2$. We split the integral
		\begin{align*}
		\int_X f(z)\frac{1}{\dist(x,z)^{n-1}}\frac{1}{\dist(y,z)^{n-1}} \di \meas(z)= I+II+III,
		\end{align*}
		where
		\begin{equation*}
			I:=\int_{B(x,r)} f(z)\frac{1}{\dist(x,z)^{n-1}}\frac{1}{\dist(y,z)^{n-1}} \di \meas(z),
		\end{equation*}
		
		\begin{equation}
			II:=\int_{B(y,r)} f(z)\frac{1}{\dist(x,z)^{n-1}}\frac{1}{\dist(y,z)^{n-1}} \di \meas(z),
		\end{equation}
		and
		\begin{equation*}
			III:=\int_{\{\dist(x,z)\geq r,\ \dist(y,z)\geq r\}} f(z)\frac{1}{\dist(x,z)^{n-1}}\frac{1}{\dist(y,z)^{n-1}} \di \meas(z).
		\end{equation*}
		  In order to estimate $I$ we observe that $\dist(y,z)\geq \dist(y,x)-\dist(x,z)\geq \dist(x,y)/2$ for all $z\in B(x,r)$, thus
		\begin{align*}
		   \int_{B(x,r)} f(z)\frac{1}{\dist(x,z)^{n-1}}\frac{1}{\dist(y,z)^{n-1}} \di   \meas(z)\leq & \frac{2^{n-1}}{\dist(x,y)^{n-1}}\int_{B(x,r)} f(z) \frac{1}{\dist(x,z)^{n-1}}\di \meas(z).
		\end{align*}
		Arguing exactly as in the proof of \eqref{a1}, we find
		$$
	    	\int_{B(x,r)} f(z) \frac{1}{\dist(x,z)^{n-1}}\di \meas(z)\leq C \dist(x,y)Mf(x), 
		$$
		and we conclude
		\begin{equation}\label{eq:estimateI}
	    	I\leq \frac{C}{\dist(x,y)^{n-2}} Mf(x).
		\end{equation}
		The estimate
		\begin{equation}\label{eq:estimateII}
	    	II\leq \frac{C}{\dist(x,y)^{n-2}} Mf(y),
		\end{equation}
		follows by the same reasoning.
		
		We are left with the estimate of $III$. By Young's inequality we have
		\begin{align*}
	    	III\leq \frac{1}{2} \int_{\{\dist(x,z)\geq r\}} f(z)\frac{1}{\dist(x,z)^{2n-2}}\di \meas(z) +\frac{1}{2} \int_{\{\dist(y,z)\geq r\}} f(z)\frac{1}{\dist(y,z)^{2n-2}} \di \meas(z).
		\end{align*}
		Using the Ahlfors regularity \eqref{eq:Ahlforsregularity}, for every $w\in X$, we obtain
		\begin{align*}
	    	\int_{\dist(w,z)\geq r} f(z)\frac{1}{\dist(w,z)^{2n-2} } &\di \meas(z) = 
	    	\sum_{k=0}^{\log_2(D/2r)} \int_{B(w,r2^{k+1})\setminus B(w,r2^k)} \frac{f(z)}{\dist(w,z)^{2n-2}} \di \meas(z)\\
	     	\leq & \sum_{k=0}^{\log_2(D/2r)}\frac{\meas(B(w,r2^{k+1}))}{(2^kr)^{2n-2}}\frac{1}{\meas(B(w,r2^{k+1}))} \int_{B(w,r2^{k+1})} f(z) \di \meas(z)\\
	    	\leq & C\sum_{k=0}^{\log_2(D/2r)}\frac{(r2^{k+1})^n}{(2^kr)^{2n-2}} Mf(w)\\
	    	\leq & \frac{C}{\dist(w,y)^{n-2}}Mf(w) \sum_{k=0}^{\log_2(D/2r)} 2^{-k(n-2)}.
		\end{align*}
		Putting this last estimate, applied with $w=x$ and $w=y$, together with \eqref{eq:estimateI} and \eqref{eq:estimateII} we obtain the desired conclusion.
	\end{proof}
	
	\begin{proof}[Proof of \autoref{lemma: key estimate}]
		First of all we remark that $|b\cdot \nabla G_x(y)+b \cdot\nabla G_y(x)|$ is well defined $\meas\times\meas$-a.e., since $b$ is a bounded  vector field and $G_x, G_y\in W^{1,p}$ for some $p>1$.
		As a first step we prove the following 
		
		\textbf{Claim:} for every $\epsilon>0$ it holds that
		\begin{equation*}
		|P_{\epsilon} (b\cdot \nabla G_x^{\epsilon})(y)+P_{\epsilon}(b \cdot \nabla G_y^{\epsilon})(x)|\leq
		\frac{C}{\dist(x,y)^{n-2}} [M(|\nabla_{\text{sym}}b|+|\div b|)(x)+M(|\nabla_{\text{sym}}b|+|\div b|)(y)],
		\end{equation*}
		for every $x,y \in X$.
		
		Recalling the result of \autoref{action of G eps} we have
		\begin{align*}
		|P_{\epsilon} (b\cdot \nabla G_x^{\epsilon})(y)+ &P_{\epsilon}(b\cdot\nabla G_y^{\epsilon})(x)|\\
		=& \left| \int_X b\cdot \nabla G_x^{\epsilon}(z)p_{\epsilon}(y,z) \di \meas(z)+\int_X b\cdot \nabla G_y^{\epsilon}(z)p_{\epsilon}(x,z) \di \meas(z) \right|\\
		=& \left| -\int_X \left[ b\cdot \nabla G_x^{\epsilon}\Delta G^{\epsilon}_y(z)+b\cdot \nabla G_y^{\epsilon}(z)\Delta G^{\epsilon}_x+\div b(G^{\epsilon}_x+G^{\epsilon}_y)\right] \di \meas(z) \right|\\
		\leq &  \left| -\int_X \left[ b\cdot \nabla G_x^{\epsilon}\Delta G^{\epsilon}_y(z)+b\cdot \nabla G_y^{\epsilon}(z)\Delta G^{\epsilon}_x-\div b(\nabla G^{\epsilon}_x\cdot \nabla G^{\epsilon}_y)\right] \di \meas(z) \right|\\
		& +\left| \int_X \div b(G^{\epsilon}_x+G^{\epsilon}_y+\nabla G^{\epsilon}_x\cdot \nabla G^{\epsilon}_y) \di \meas(z)\right|\\
		=&2\left| \int_X \nabla_{\text{sym}}b(\nabla G^{\epsilon}_x, \nabla G^{\epsilon}_y) \di \meas(z)\right| + \left| \int_X \div b(G^{\epsilon}_x+G^{\epsilon}_y+\nabla G^{\epsilon}_x\cdot \nabla G^{\epsilon}_y) \di \meas(z)\right|.
		\end{align*}
		Now using \eqref{G-d estimate 1} and \eqref{ nabla G estimate} we get
		\begin{align*}
		|P_{\epsilon} (b\cdot \nabla G_x^{\epsilon})(y)+ &P_{\epsilon}(b\cdot \nabla G_y^{\epsilon})(x)|\\
		\le&C\left| \int_X(|\nabla_{\text{sym}}b|(z)+|\div b(z)|)\frac{1}{d(x,z)^{n-1}}\frac{1}{\dist(y,z)^{n-1}} \di \meas(z)\right|,
		\end{align*}
		where we have implicitly exploited the inequality  $\frac{1}{\dist(\cdot,z)^{n-2}}\leq\frac{D}{\dist(\cdot,z)^{n-1}}$.
		Applying \autoref{maximal estimate} with $f:=|\nabla_{\text{sym}}b|+|\div b|$ we conclude the proof of the claim.
		
		We want to prove now that 
		$$
		P_{\epsilon} (b\cdot \nabla G_x^{\epsilon})(y)+P_{\epsilon}(b\cdot \nabla G_y^{\epsilon})(x)\rightarrow b\cdot \nabla G_x(y)+b\nabla \cdot G_y(x)
		$$
		strongly in $L^p(X\times X, \meas\times \meas)$ when $\epsilon\to 0$; this convergence result together with the uniform estimate we proved above will yield the desired conclusion (by considering a sequence $(\epsilon_i)_{i\in\setN}$ such that $\epsilon_i\to 0$ and the above considered convergence holds in the $\meas\times\meas$-a.e. sense).
		% %TODO:sottolineato anche qua sopra che bisogna considerare una sottosuccessione che converge quasi ovunque. 
		
		Entering into the details we are going to prove that 
		$$
		\int_X\int_X |P_{\epsilon} (b\cdot \nabla G_x^{\epsilon})(y)-b\cdot \nabla G_x(y)|^p \di \meas(x)\di \meas(y)
		\rightarrow 0\quad   \left(\epsilon\to 0\right).    
		$$

		Recalling the $L^p$-norm contractivity property of the semigroup $P_t$ we have that for any fixed $x\in X$ it holds 
		\begin{align*}
		\norm{P_{\epsilon} (b\cdot \nabla G_x^{\epsilon})-b\cdot \nabla G_x}_{L^p}
		\leq &
		\norm{P_{\epsilon} (b\cdot \nabla G_x^{\epsilon})-P_{\epsilon}(b\cdot \nabla G_x)}_{L^p}+\norm{P_{\epsilon} (b\cdot \nabla G_x)-b\cdot \nabla G_x}_{L^p}\\
		\leq &
		 \norm{\abs{b}}_{L^{\infty}}\norm{\nabla (G_x^{\epsilon}-G_x)}_{L^p}+\norm{P_{\epsilon} (b\cdot \nabla G_x)-b\cdot \nabla G_x}_{L^p}.
		\end{align*}
		The last two terms go to zero when $\epsilon\to 0$, moreover they are uniformly bounded in $x$, thus
		
		\begin{align*}
		\int_X\int_X |P_{\epsilon} (b\cdot \nabla    &G_x^{\epsilon})(y)-b\cdot \nabla G_x(y)|^p \di \meas(y)\di \meas(x) \\
		\leq &\int_X \norm{b}_{L^{\infty}}\norm{\nabla (G_x^{\epsilon}-G_x)}_{L^p}\di \meas(x)+\int_X \norm{P_{\epsilon} (b\cdot \nabla G_x)-b\cdot \nabla G_x}_{L^p} \di \meas(x)\\
		\end{align*}
		goes to zero by the dominated convergence theorem.
	\end{proof}

	\subsection{A Lusin-type regularity result}\label{def:RLF}
	
	Throughout this section the time dependent vector field $b_t$ and the regular Lagrangian flow $\XX_t$ associated to $b_t$, with compressibility constant $L$, are fixed. Our aim is to implement a strategy very similar to the one developed in \cite{CrippaDeLellis08} in order to prove our main regularity result \autoref{thm: Lusin-type regularity of RLF}.
	
	We begin by observing that the results of \autoref{prop: properties of G} ensure that, possibly increasing the constant $A$ and setting $\bar{G}(x,y):=G(x,y)+\bar{A}$ we have
	\begin{equation}\label{eq:estimatebarG}
	\frac{1}{A\dist(x,y)^{n-2}}\le \bar{G}(x,y)\le \frac{A}{\dist(x,y)^{n-2}},
	\end{equation}
	for any $x,y\in X$ such that $x\neq y$. It follows in particular that $\bar{G}(x,y)>\alpha>0$ for any $x,y\in X$ such that $x\neq y$.

	Observe that, in terms of the function $\bar{G}$, the statement of \autoref{lemma: key estimate} can be rewritten as 
	\begin{equation}\label{eq:maximalestimateGbar}
	|b\cdot \nabla \bar{G}_x(y)+b \cdot \nabla \bar{G}_y(x)|\leq
	C\bar{G}(x,y) [M(|\nabla_{\text{sym}}b|+|\div b|)(x)+M(|\nabla_{\text{sym}}b|+|\div b|)(y)],
	\end{equation}
	for $\meas\times\meas$-a.e. $(x,y)\in X\times X$.

    We introduce, for any $t\in[0,T]$ and for any $0<r\le D$, the functional
	\begin{equation*}
	Q_{t,r}(x):=\dashint_{B(x,r)} \log \left( 1+\frac{1}{A}\left(\frac{ \dist(\XX_t(x),\XX_t(y))}{r}\right)^{n-2} \right) \di \meas(y),
	\end{equation*}
	where $A$ is the constant introduced in \eqref{eq:estimatebarG}.	
	Moreover we set
	\begin{equation}\label{def: Q star}
	Q^*(x):=\sup_{0\leq t\leq T}\sup_{0<r\leq D}Q_{t,r}(x).
	\end{equation}
	
	With the aim of finding bounds on $Q^*$, we first state and prove a technical lemma.
	
	\begin{lemma}
		\label{lemma: faticaccia} Assume that $b\in L^1((0,T), L^2(TM))$ and that it is bounded. Then, for $\meas\times \meas$-a.e. $(x,y)\in X\times Y$ the map	$t\to G(\XX_t(x),\XX_t(y))$ belongs to $W^{1,1}((0,T))$ and its derivative is given by the formula
        \begin{equation}\label{z1}
		\frac{d}{dt}  G(\XX_t(y),\XX_t(x))= b_t\cdot \nabla G_{\XX_t(x)}(\XX_t(y))+b_t\cdot \nabla G_{\XX_t(y)}(\XX_t(x)),
		\end{equation}
		for a.e. $t\in (0,T)$.
	\end{lemma}
	
	\begin{proof}
	
		From \eqref{eq:spectralidentity} and the definition of $G^{\epsilon}$ we find the pointwise identity
		\begin{equation}\label{eq:G eps eigenfunctions}
		G^{\epsilon}(x,y)=\sum_{i=0}^{\infty}\frac{e^{-\lambda_i\eps}}{\lambda_i} u_i(x)u_i(y),
		\end{equation}
		for every $x,y\in X$. Setting
		\begin{equation*}
			G^{\eps,N}(x,y):=\sum_{i=0}^N\frac{e^{- \lambda_i \eps}}{\lambda_i} u_i(x)u_i(y),
		\end{equation*}
	     we have that, for $\meas\times \meas$-a.e. $(x,y)\in X\times X$, the map $t\mapsto G^{\eps,N}(\XX_t(x),\XX_t(y))$
		is absolutely continuous for every $N\in \setN$. Moreover, for a.e. $t\in(0,T)$, it holds
		\begin{align}\label{z}
		  \frac{d}{dt}G^{\eps,N}(\XX_t(x),\XX_t(y))= &\sum_{i=0}^N\frac{e^{-\epsilon \lambda_i}}{\lambda_i}(b_t\cdot\nabla u_i(\XX_t(x))u_i(\XX_t(y))
		  +u_i(\XX_t(x))b_t\cdot\nabla u_i(\XX_t(y)))\\
		  =& b_t\cdot\nabla G^{\eps,N}_{\XX_t(x)}(\XX_t(y))+b_t\cdot\nabla G^{\eps,N}_{\XX_t(y)}(\XX_t(x)),
		\end{align}
		since $u_i\in \Test(X,\dist,\meas)$ for any $i\in\setN$. Our aim is to pass to the limit \eqref{z}, first letting $N\to \infty$ and then $\eps\to 0$.  Observe that, for every $\eps>0$, when $N\to \infty$ we have $G^{\eps,N} \to G^{\eps}$ in $W^{1,2}(X\times X)$ (moreover it holds $G_x^{\eps,N}\to G_x$ strongly in $W^{1,2}(X)$, uniformly in $x\in X$) and also, when $\eps\to 0$, $G^{\eps}\to G$ in $L^1(X\times X)$, and $G_x^{\eps}\to G_x$ in $W^{1,p}(X)$ for every $p\in [1,\frac{n}{n-1})$, uniformly in $x\in X$ (see \autoref{cor: G sobolev}).
	    Moreover $G^{\eps}\in \Test(X\times X)$ (it can be proved arguing as in \autoref{action of G eps}).
	    
	    With this said, in order to conclude the proof, it suffices to show the following technical result: let $(F^n(x,y))_{n\in \mathbb{N}}$ be a sequence of symmetric functions belonging to $\Test(X\times X)$, assume that $F^n$ satisfies \eqref{z} for every $n\in \mathbb{N}$. If $F^n\to F$ in $L^1(X\times X)$ and there exists $p>1$ such that, for every $x\in X$, $F^n(x,\cdot):=F_x^n(\cdot)\to F_x(\cdot)$ in $W^{1,p}(X)$, uniformly w.r.t. $x\in X$, then $t\to F(\XX_t(x),\XX_t(y))$ belongs to $W^{1,1}((0,T))$ for $\meas\times \meas$-a.e. $(x,y)\in X\times X$ and satisfies \eqref{z}. 
	    
	    Let us fix $t\in [0,T]$, starting from the $\meas\times \meas$-a.e. equality
	    \begin{equation}
	    	F^n(\XX_t(x),\XX_t(y))-F^n(x,y)=\int_0^t \left\lbrace b_s\cdot \nabla F^n_{\XX_s(x)}(\XX_s(y))+b_s\cdot \nabla F^n_{\XX_s(y)}(\XX_s(x))) \right\rbrace \di s, 
	    \end{equation}
	    we wish to pass to the limit for $n\to \infty$. Observe that the left hand side converges to $F(\XX_t(x),\XX_t(y))-F(x,y)$ in $L^1(X\times X)$ (here the compressibility property of $\XX_t$ plays a role), it remains only to prove that the right hand side converges to 
	    \[
	    \int_0^t \left\lbrace b_s\cdot \nabla F_{\XX_s(x)}(\XX_s(y))+b_s\cdot \nabla F_{\XX_s(x)}(\XX_s(y)) \right\rbrace \di s,
	    	\qquad \text{in}\ L^1(X\times X).
	    \]
	    Using again the compressibility property of $\XX_t$ we have
	    \begin{align*}
	    \int_{X\times X} & \abs{\int_0^t b_s\cdot \nabla F^n_{\XX_s(x)}(\XX_s(y))\di s-\int_0^t b_s\cdot \nabla F_{\XX_s(x)}(\XX_s(y))\di s}^p  \di \meas(x)\di \meas(y)\\
	    \leq & t^{1-1/p}\norm{b}_{L^{\infty}} \int_0^t\int_X \int_X |\nabla ( F^n_{\XX_s(x)}- F_{\XX_s(x)})|^p(\XX_s(y)) \di \meas(y)\di \meas(x)\di s\\
	    \leq & L^2t^{2-1/p} \norm{b}_{L^{\infty}} \int_X\norm{\nabla (F^n_x-F_x)}_{L^p}^p \di \meas(x),
	    \end{align*}
	    that, under our assumptions on the sequence $F^n$, goes to zero when $n\to \infty$. Arguing similarly for the term $\int_0^t\nabla b_s \cdot F^n_{\XX_s(y)}(\XX_s(x)) \di s$ we conclude the proof.
	\end{proof}

	\begin{theorem}\label{th 1}
		Let $b$ be a time dependent vector field, assume that $\abs{\nabla_{\text{sym}}b}$ and $\div b$ belong to $L^1((0,T);L^2(X,\meas))$. Let $\XX_t$ be a Regular Lagrangian flow associated to $b$, with compressibility constant $L$.
		Then, with the above introduced notation, we have
		\begin{equation*}
		\norm{Q^*}_{L^2}\leq C \left[ L \int_0^T \norm{|\nabla_{\text{sym}}b_s|+|\div   b_s|}_{L^2} \di s+1 \right],
		\end{equation*}
		where $C=C(T,X,\dist,\meas)$.
	\end{theorem}

	\begin{proof}
		
		As a first step we introduce the functional
		\begin{equation}\label{CDL functionals}
		    \Phi_{t,r}(x):=\dashint_{B(x,r)}\log \left(   1+\frac{1}{r^{n-2}\bar{G}(\XX_t(x),\XX_t(y)))} \right) \di \meas(y),
		\end{equation}
		for $r\in (0,D)$ and $t\in [0,T]$.
		
		We observe that the monotonicity of the logarithm and our construction grant that $Q_{t,r}\leq\Phi_{t,r}$ pointwise for any $t\in[0,T]$ and for any $0<r\le D$.
		
		What we just proved tells that it suffices to bound $\Phi_{t,r}$ (which in some sense is a more ``regular'' functional) in order to bound $Q_{t,r}$. To this aim we fix $r>0$ and $t\in [0,T]$. 
		By \autoref{lemma: faticaccia} we have that $t\to \Phi_{t,r}(x)$ belongs to $W^{1,1}((0,T))$ for $\meas$-a.e. $x\in X$ (actually it is absolutely continuous since it is continuous) and it holds
		\begin{align*}
		\Phi_{t,r}(x)=&\Phi_{0,r}(x)+ \int_0^t \frac{d}{ds}\Phi_{s,r}(x) \di s\\
		\leq &\Phi_{0,r}(x)+\int_0^t\dashint_{B(x,r)} \frac{|\frac{d}{ds}\bar{G}(\XX_s(x),\XX_s(y))|}{\bar{G}(\XX_s(x),\XX_s(y))}\cdot
		\frac{1}{\bar{G}(\XX_s(x),\XX_s(y))r^{n-2}+1} \di \meas(y)\di s\\
		\leq & \Phi_{0,r}(x)+\int_0^t\dashint_{B(x,r)} \frac{| b_s\cdot\nabla \bar{G}_{\XX_s(x)}(\XX_s(y))+b_s\cdot\nabla \bar{G}_{\XX_s(y)}(\XX_s(x))|}{\bar{G}(\XX_s(x),\XX_s(y))} \di \meas(y)\di s.
		\end{align*}	
		Setting $g_s:=|\nabla_{\text{sym}}b_s|+|\div b_s|$, from \autoref{lemma: key estimate} and \eqref{eq:maximalestimateGbar} we obtain
		\begin{align*}
		\Phi_{t,r}(x)\leq & \int_0^t \dashint_{B(x,r)} \frac{| b_s\cdot \nabla \bar{G}_{\XX_s(x)}(\XX_s(y))+b_s\cdot \nabla \bar{G}_{\XX_s(y)}(\XX_s(x))|}{\bar{G}(\XX_s(x),\XX_s(y))} \di \meas(y) \di s+\Phi_{0,r}\\
		\leq & C\int_0^t \dashint_{B(x,r)} \left(Mg_s(\XX_s(x))+Mg_s(\XX_s(y))\right) \di \meas(y) \di s+\Phi_{0,r}\\
		\leq & C\int_0^t [Mg_s(\XX_s(x))+M(Mg_s\XX_s(\cdot))(x)] \di s+1/A
		\end{align*}
		for $\meas$ a.e. $x\in X$, and the negligible set does not depend on $r$. 
		From \eqref{G-d estimate 2} we get 
		\begin{align*}
		\sup_{0\leq t\leq T}\sup_{0<r\leq D}Q_{t,r}(x)
		& \leq \sup_{0\leq t\leq T}\sup_{0<r\leq D}\Phi_{t,r}(x)\\
		&\leq C\int_0^T [Mg_s(\XX_s(x))+M(Mg_s(\XX_s(\cdot))(x)] \di s+1/A,
		\end{align*}
		for $\meas$ a.e. $x\in X$. Taking the $L^2$-norms (here the assumption that the RLF has compressibility constant $L<+\infty$ enters into play once more) we obtain our thesis.	  
	\end{proof}

	Below we state and prove the main regularity result for regular Lagrangian flows.
	
	\begin{theorem}\label{thm: sobolev vectorfield, regularity}
		Let $b$ be a time dependent vector field and let $\XX$ be a regular Lagrangian flow associated to $b$ with compressibility constant $L$. For any $x,y\in X$ and for any $t\in[0,T]$ it holds
		\begin{equation}\label{lip regularity of X 1}
		\dist(\XX_t(x),\XX_t(y))\le C e^{C\left(Q^*(x)+Q^*(y)\right)}\dist(x,y),
		\end{equation}
		where $Q^*$ has been defined in \eqref{def: Q star}, and $C=C(X,\dist,\meas)$.
		
		Moreover, when $b$ is bounded and $\abs{\nabla_{\text{sym}}} b, \div b \in L^1((0,T); L^2(X,\meas))$, for every $\epsilon>0$ there exists a Borel set $E\subset X$ such that
		$\meas(X\setminus E)<\epsilon$ and for every $x,y\in E$
		\begin{equation}\label{eq:Lipregularity}
		\dist(\XX_t(x),\XX_t(y))\leq C\exp\left(2C\frac{\norm{Q^*}_{L^2}}{\sqrt{\epsilon}}\right) \dist(x,y)
		\end{equation}
		for every $t\in [0,T]$, where we remark that this last statement makes sense since, under our regularity assumptions on $b$, \autoref{th 1} grants that $\norm{Q^*}_{L^2}<+\infty$. 
	\end{theorem}
	
	\begin{proof}
		Fix any $x,y\in X$ such that $x\neq y$ and set $r:=\dist(x,y)$. Exploiting the inequality $(a+b)^m\le 2^{m-1}(a^m+b^m)\le 2^m(a^m+b^m)$, the triangular inequality and the subadditivity and monotonicity of $t\mapsto\log(1+t)$, we obtain that
		\begin{align*}
		\log \left( 1+\frac{1}{A}\left(\frac{ \dist(\XX_t(x),\XX_t(y))}{2r}\right)^{n-2} \right) 
		\leq &\log \left( 1+\frac{1}{A}\left(\frac{ \dist(\XX_t(x),\XX_t(z))  }{r}\right)^{n-2} \right)\\
		&+\log \left( 1+\frac{1}{A}\left(\frac{ \dist(\XX_t(z),\XX_t(y))}{r}\right)^{n-2} \right),
		\end{align*}
		for any $z\in X$. Let us fix $w\in X$ such that $\dist(x,w)=\dist(y,w)=r/2$ (observe that such a point exists, since $(X,\dist)$ is a geodesic metric space) and we take the mean value of the above written inequality (w.r.t. the $z$ variable) over $B(w,r/2)$ obtaining
		\begin{align*}
		\log \left( 1+\frac{1}{A}\left(\frac{ \dist(\XX_t(x),\XX_t(y))}{2r}\right)^{n-2} \right)
		\le&	\dashint_{B(w,r/2)}\log \left( 1+\frac{1}{A}\left(\frac{ \dist(\XX_t(x),\XX_t(z))  }{r}\right)^{n-2} \right)\di\meas(z)\\
		&+\dashint_{B(w,r/2)}\log \left( 1+\frac{1}{A}\left(\frac{ \dist(\XX_t(z),\XX_t(y))}{r}\right)^{n-2} \right)\di\meas(z)\\
		\le& C\dashint_{B(x,r)}\log \left( 1+\frac{1}{A}\left(\frac{ \dist(\XX_t(x),\XX_t(z))  }{r}\right)^{n-2} \right)\di\meas(z)\\
		&+ C\dashint_{B(y,r)}\log \left( 1+\frac{1}{A}\left(\frac{ \dist(\XX_t(z),\XX_t(y))}{r}\right)^{n-2} \right)\di\meas(z),
		\end{align*}
		where in the second inequality we enlarge the domain of integration and control the ratios between volumes of balls with radii $r/2$ and $r$ in a uniform way, thanks to the Ahlfors regularity assumption.	
		
		It follows by the definition of $Q^*$ that for any $x,y\in X$ such that $x\neq y$ and for any $t\in[0,T]$ it holds
		\begin{equation*}
		\log\left(1+\frac{1}{A}\left(\frac{\dist(\XX_t(x),\XX_t(y))}{2\dist(x,y)}\right)^{n-2}\right)\le C\left(Q^*(x)+Q^*(y)\right),
		\end{equation*} 			
		which easily yields \eqref{lip regularity of X 1}. 		 
		
		Let us define $E:=\set{x\in X\ :\ Q^*(x)\leq \norm{Q^*}_{L^2}/\sqrt{\epsilon}}$, by Chebyshev inequality we deduce that $\meas(X\setminus E)<\epsilon$. 
		The conclusion of \eqref{eq:Lipregularity} now directly follows from \eqref{lip regularity of X 1}.
	\end{proof}

     \subsection{The case $n = 2$}\label{sec: n leq 2}

     In order to conclude the proof of our result \autoref{thm: Lusin-type regularity of RLF} we have to deal with the case $n=2$.
     
     In order to reduce this case to an application of the result we proved for $n>2$ we ``add a dimension'' to the given space by considering its product with the standard $\mathbb{S}^1$. 
          
     To this aim, given a $2$-Ahlfors regular $\RCD^*(K,N)$ m.m.s. $(X,\dist,\meas)$ we define $(\bar{X},\bar{\dist},\bar{\meas})$ by
      \begin{itemize}
     	\item[1)] $\bar{X}:=X\times \mathbb{S}^1$;
     	\item[2)] $\bar{\dist}^2((x,s), (x',s')):=\dist^2_{X}(x,x')+\dist^2_{\mathbb{S}^1}(s,s')$ for every $x,x'\in X$ and $s,s'\in \mathbb{S}^1$;
     	\item[3)] $\bar{\meas}:= \meas \times \di s$, where $\di s$ is the (normalized) volume measure of $\mathbb{S}^1$.
     \end{itemize}
     Then $(\bar{X},\bar{\dist},,\bar{\meas})$ is an $\RCD^*(K,N+1)$ m.m.s. (see \cite[sect. 6]{AmbrosioGigliSavare14} and \cite{BacherSturm10}) and it is $3$-Ahlfors regular, as an elementary application of Fubini's theorem shows.
     
     We will denote by $\pi_1$ and $\pi_2$ the canonical projections from $\bar{X}$ to $X$ and $\mathbb{S}^1$ respectively. With this said we introduce the so-called algebra of tensor products by
     \begin{equation*}
     \mathcal{A}:=\left\lbrace \sum_{j=1}^ng_j\circ\pi_1h_j\circ\pi_2: n\in\setN, g_j\in W^{1,2}\cap L^{\infty}(X)\text{ and }h_j\in W^{1,2}\cap L^{\infty}(\mathbb{S}^1) \forall j=1,\dots,n\right\rbrace. 
     \end{equation*}
     
     A crucial property for the rest of the discussion in this section is the strong form of density of the product algebra $\mathcal{A}$ (the terminology is borrowed from \cite{GigliRigoni17}), namely it holds that for any $f\in W^{1,2}(\bar{X})\cap L^{\infty}(X)$ there exists a sequence $(f_n)_{n\in\setN}$ with $f_n\in\mathcal{A}$ uniformly bounded and converging to $f$ in $W^{1,2}(\bar{X})$. This property can be proved with minor modifications of the strategy developed in \cite{GigliHan15}, where the case of products with the Euclidean line or intervals is considered.
     We also remark that a more direct approach to the proof of this density result can be obtained exploiting the result of \cite[Theorem B.1]{AmbrosioStraTrevisan17}. Indeed, knowing that the algebra generated by distances from points is dense, the observation that the distance squared from a point in the product belongs to the product algebra (actually $\bar{\dist}^2((x,s),\cdot)=\dist_X^2(x,\cdot)+\dist_{\mathbb{S}^1}^2(s,\cdot)$) together with an approximation procedure (which is needed to recover the distance from the distance squared) yields the desired conclusion.

     In order to be able to apply \autoref{thm: Lusin-type regularity of RLF} in the space $\bar{X}$ we are going to lift the given vector field $b_t$ on $X$ to a vector field $\bar{b}_t$ on $\bar{X}$. 
     
     In \cite{GigliRigoni17} the study of tangent and cotangent moduli of product spaces is performed in great generality. We just observe here that, in the case of our interest, we are in position to lift $b_t$ in a trivial way by saying that, for any $f\in W^{1,2}(\bar{X})$,
     \begin{equation*}
     \bar{b}_t\cdot\nabla f(x,s)=b_t\cdot\nabla f_s(x) 
     \end{equation*}
     for $\bar{\meas}$-a.e. $(x,s)\in\bar{X}$, where $f_s(x):=f(s,x)$ and we are implicitly exploiting the tensorization property of the Cheeger energy (see \cite{AmbrosioGigliSavare14}). Observe that if $b_t\in L^{\infty}((0,T); L^{\infty}(TX))$ then $\bar{b}_t\in L^{\infty}((0,T);L^{\infty}(T\bar{X}))$ and the norms are actually preserved.

     % %Ho aggiunto un commento sulla norma del campo esteso, che viene preservata.
     % %bisognerebbe anche osservare che si mantiene la misurabilità rispetto al tempo...

     % %Cambio un po' la presentazione, non direi che definiamo un RLF, direi che definiamo una famiglia di mappe e che verifichiamo che si tratta di un RLF. Poi isolo questo fatto e il fatto che il campo liftato ha espressioni esplicite per derivata simmetrica e divergenza in due lemmi, la divisione in step della discussione mi sembra poco comune nei lavori che si vedono in giro.
     
     Given a RLF $\XX_t$ of $b_t$ on $(X,\dist,\meas)$ we go on by setting
     \begin{equation*}
     	   \bar{\XX}_t(x,s):=(\XX_t(x),s) \qquad\qquad \forall x\in X,\ \forall s\in \mathbb{S}^1,
     \end{equation*}
     for every $t\in [0,T]$.
     
    In \autoref{lemma:liftedisRLF} and \autoref{lemma:liftedvfisregular} below we prove that $\XX_t$ is a Regular Lagrangian flow of $\bar{b}_t$ and that $\bar{b}_t$ inherits the Sobolev regularity from $b_t$. These remarks will put us in position to apply \autoref{thm: sobolev vectorfield, regularity}.
    
     \begin{lemma}\label{lemma:liftedisRLF}
     Under the previous assumptions $(\bar{\XX}_t)_{t\in[0,T]}$ is a Regular Lagrangian flow of $\bar{b}_t$.
     \end{lemma}
     
     \begin{proof}
     We begin by observing that, if $\XX_t$ has compression bounded by $L$, then $\bar{\XX}_t$ has compression bounded by $L$ itself, as a simple application of the change of variables formula shows.
     
     With this said it remains to check that condition 3 in \autoref{def:Regularlagrangianflow} is satisfied (since the trajectories of $\bar{\XX_t}$ inherit the continuity from the trajectories of $\XX_t$). To this aim let $\tilde{\mathcal{A}}$ be the algebra generated by functions of the form $f\circ\pi_1g\circ\pi_2$, where $g\in\Test(X,\dist,\meas)$ and $f\in\Test(\mathbb{S}^1)$.
     Observe that, for $f\in\tilde{\mathcal{A}}$, the identity
     \begin{equation}\label{eq:checkRLF}
     \frac{\di}{\di t} f(\bar{\XX}_t(x,s))=\bar{b}_t\cdot \nabla f (	\bar{\XX}_t(x,s))
     \end{equation}
     for $\bar{\meas}$-a.e. $(x,s)$ and for $\leb^1$-a.e. $t\in(0,T)$
     directly follows from the assumption that $\XX_t$ is a Regular Lagrangian flow of $b_t$ on $X$ and from the definition of $\bar{b}_t$.
     To conclude that \eqref{eq:checkRLF} holds true for any $f\in \Test(\bar{X})$ it is now sufficient to take into account the density of the product algebra in the strong form. 
     \end{proof}

     \begin{lemma}\label{lemma:liftedvfisregular}
     Assume that $b\in L^2(TX)$ has divergence in $L^2(\meas)$. Then $\bar{b}$ has divergence in $L^2(\bar{\meas})$ and it holds $\div\bar{b}(x,s)=\div b(x)$ for $\bar{\meas}$-a.e. $(x,s)\in \bar{X}$.
     
     Moreover, if $b$ has symmetric derivative in $L^2(\meas)$ according to \autoref{def: symder Elia-Dan}, then $\bar{b}$ has symmetric derivative in $L^2(\bar{\meas})$ itself and it holds $\norm{\nabla_{sym}\bar{b}}_{L^2(\bar{\meas})}\le\norm{\nabla_{sym}b}_{L^2(\meas)}$, where we omitted the implicit dependence from the space of the divergence and the symmetric covariant derivative.
     \end{lemma}
     
     \begin{proof}
     The proof of the first conclusion can be found in \cite[Proposition 3.15]{GigliRigoni17}. 
     
     We pass to the proof of the second statement, which is strongly inspired by the proof of an analogous result concerning the Hessian that can be found in \cite[Appendix A]{GigliRigoni17}.
     
     Recall that we defined $\abs{\nabla_{sym}\bar{b}}$ to be the function $h\in L^{2}(\bar{\meas})$ with smallest $L^2$-norm such that 
     \begin{equation}\label{eq:symderoflifted}
       \frac{1}{2}\left|\int_{\bar{X}} \bar{b}\cdot \nabla g \Delta f + \bar{b}\cdot \nabla f \Delta g- \div \bar{b} (\nabla f\cdot \nabla g) \di \bar{\meas}\right|  \leq \int_{\bar{X}} h |\nabla f||\nabla g| \di \bar{\meas}
     \end{equation}
   for any $f,g\in\Test(\bar{X})$. Therefore to prove the desired conclusion it suffices to show that $h(x,s):=\abs{\nabla_{sym}b}(x)$ is an admissible function in \eqref{eq:symderoflifted}. Moreover, thanks to the strong density of the algebra $\mathcal{A}$ and to the approximation result of \cite[Lemma A.3]{GigliRigoni17}, it is sufficient to verify \eqref{eq:symderoflifted} in the case where $f,g\in\tilde{\mathcal{A}}$, where the algebra $\tilde{\mathcal{A}}$ was introduced in the proof of \autoref{lemma:liftedisRLF}.
   
   Denoting by $\Delta_X$ and $\Delta_{\mathbb{S}^1}$ the Laplacians on $X$ and $\mathbb{S}^1$ respectively, we recall from \cite[pg. 52]{AmbrosioGigliSavare14} that (with a slight abuse of notation) it holds $\Delta_{\bar{X}}=\Delta_X+\Delta_{\mathbb{S}^1}$.
   
   Then we compute  
   \begin{align*}
       	\frac{1}{2}\int_{\bar{X}}  \bar{b} \cdot \nabla f \Delta g+\bar{b} &\cdot \nabla g \Delta f-\div \bar{b} \nabla f\cdot \nabla g   \di \bar{\meas}\\
       	=&\frac{1}{2}\int_{\mathbb{S}^1} \left[ \int_X b\cdot \nabla_X f \Delta_X g+b\cdot \nabla_X g \Delta_X f-\div b\ \nabla_X f\cdot \nabla_X g \di \meas \right] \di s\\
       	&+ \frac{1}{2}\int_{\mathbb{S}^1} \left[ \int_X b\cdot \nabla_X f \Delta_{\mathbb{S}^1} g+b\cdot \nabla_X g \Delta_{\mathbb{S}^1} f-\div b\ \nabla_{\mathbb{S}^1} f\cdot \nabla_{\mathbb{S}^1} g \di \meas\right] \di s,
       \end{align*}
   where we exploited the definition of $\bar{b}$, the previously proven identity $\div\bar{b}=\div b\circ\pi_1$ and the tensorization of the Cheeger energy again. The first of the two terms appearing above is bounded by 
   \begin{equation*}
   \int_{\bar{X}}h\abs{\nabla_X f}\abs{\nabla_Y g}\di\bar{\meas}\le\int_{\bar{X}}h\abs{\nabla f}\abs{\nabla g}\di\bar{\meas},
   \end{equation*}
   since $b$ has symmetric derivative in $L^2$.
   
   To conclude we are going to prove that 
   \begin{equation*}
      R:=\int_{\mathbb{S}^1} \left[ \int_X b\cdot \nabla_X f \Delta_{\mathbb{S}^1} g+b\cdot \nabla_X g \Delta_{\mathbb{S}^1} f-\div b\ \nabla_{\mathbb{S}^1} f\cdot \nabla_{\mathbb{S}^1} g \right] \di \meas \di s=0.
   \end{equation*}  
     
   Observe that applying the Leibniz rule for the divergence and integrating by parts we obtain
     	 \begin{align*}
     	 \int_{\mathbb{S}^1} \int_X b\cdot \nabla_X f \Delta_{\mathbb{S}^1} g \di \meas \di s
         	 =& -\int_{\mathbb{S}^1}\int_X \div b\  f\ \Delta_{\mathbb{S}^1} g \di\meas \di s 
         	 +\int_{\mathbb{S}^1}\int_X \div(bf)\Delta_{\mathbb{S}^1}g \di \meas \di s\\
         	 =&\int_X \div b\  \int_{\mathbb{S}^1} \nabla_{\mathbb{S}^1}f \cdot\nabla_{\mathbb{S}^1} g \di s\di \meas
         	 -\int_{\mathbb{S}^1}\int_X f\  b\cdot \nabla_X (\Delta_{\mathbb{S}^1} g) \di \meas \di s\\
         	 =& \int_{\mathbb{S}^1}\int_X \div b\ \nabla_{\mathbb{S}^1}f\cdot \nabla_{\mathbb{S}^1} g \di \meas \di s
         	 -\int_{\mathbb{S}^1}\int_X f\  b\cdot \nabla_X (\Delta_{\mathbb{S}^1} g) \di \meas \di s,
         \end{align*}
     which yields to
     
     \begin{equation*}
        	    R=\int_{\mathbb{S}^1}\int_X f \left[\Delta_{\mathbb{S}^1}(b\cdot \nabla_X g)-b\cdot \nabla_X (\Delta_{\mathbb{S}^1} g)\right] \di \meas \di s.
        \end{equation*}
   To get the desired conclusion we just observe that, for any $g\in\tilde{\mathcal{A}}$, it holds
     
     \begin{equation*}
     \Delta_{\mathbb{S}^1}(b\cdot \nabla_X g)=b\cdot \nabla_X (\Delta_{\mathbb{S}^1} g).
     \end{equation*}
     
     \end{proof}
     
    % %Rispetto a quanto avevo detto non ho messo l'uguaglianza dei moduli delle derivate simmetriche perché non mi è del tutto evidente il motivo per cui l'oggetto che definiamo debba essere minimale puntualmente oltre che rispetto alla norma L^2. Questo ingrediente mi sembra necessario per ottenere l'uguaglianza (ma anche solo la disuguaglianza puntuale quasi ovunque invece che solo la disuguaglianza fra le norme). 

   As we anticipated the results of \autoref{lemma:liftedisRLF} and \autoref{lemma:liftedvfisregular} put us in position to apply \autoref{thm: Lusin-type regularity of RLF} to $\bar{\XX}_t$. 
   
   It follows that there exist a function $Q^*(x,s)$ and a constant $C=C(X,\dist,\meas)$ such that
    \begin{equation*}
        \bar{\dist}(\bar{\XX}_t(x,s),\bar{\XX}_t(x',s'))\leq Ce^{C(Q^*(x,s)+Q^*(x',s'))} \bar{\dist}((x,s),(x',s')),
    \end{equation*}
    for every $x,x'\in X$ and $s,s'\in \mathbb{S}^1$. Choosing $s=s'$ and setting $Q^*(x):=\sup_{s\in \mathbb{S}^1} Q^*(x,s)$ we obtain that
    \begin{equation*}
         \dist(\XX_t(x),\XX_t(y))\leq Ce^{C(Q^*(x)+Q^*(y))} \dist(x,y),
    \end{equation*}
    for any $x,y\in X$.
	Moreover, it follows from the proof of \autoref{th 1} and from the result (and the proof) of \autoref{lemma:liftedvfisregular} that
	\begin{equation}
      \sup_{s\in S^1} Q^*(x,s)\leq C\int_0^T[ Mg_t(\XX_t(x))+M(Mg_t(\XX_t(\cdot)))(x)] \di s+C,
	\end{equation}
	where $g:=|\nabla_{\text{sym}} b_t|(x)+|\div b_t|(x)$. Thus there exists  $C=C(T,X,\dist,\meas)$ such that
    \begin{equation*}
       \norm{Q^*}_{L^2}\leq C \left[ L \int_0^T \norm{|\nabla_{\text{sym}}b_t|+|\div   b_t|}_{L^2} \di t+1 \right].
    \end{equation*}
    The regularity result is now proved.
    
    % %qua forse bisogna dire esplicitamente che si può fare il conto con una cosa che non è il modulo della derivata simmetrica ma è ammissibile nella definizione di derivata simmetrica

	% %TODO: uniformare la bibliografia

\end{document}